\documentclass[12pt]{article}
\usepackage{amssymb,amsmath, txfonts}
\begin{document}

\def\abstractname{\bf Abstract}
\def\dfrac{\displaystyle\frac}
\let\oldsection\section
\renewcommand\section{\setcounter{equation}{0}\oldsection}
\renewcommand\thesection{\arabic{section}}
\renewcommand\theequation{\thesection.\arabic{equation}}
\newtheorem{theorem}{\indent Theorem}[section]
\newtheorem{lemma}{\indent Lemma}[section]
\newtheorem{proposition}{\indent Proposition}[section]
\newtheorem{definition}{\indent Definition}[section]
\newtheorem{remark}{\indent Remark}[section]
\newtheorem{corollary}{\indent Corollary}[section]
\def\pd#1#2{\displaystyle\frac{\partial#1}{\partial#2}}
\def\d#1{\displaystyle\frac{d#1}{dt}}

\title{\LARGE\bf Global solvability and stability
to a nutrient-taxis model with porous medium slow diffusion
\thanks{This work is supported by NSFC(11471127, 11771156, 11571380), Guangdong Natural Science Funds for
Distinguished Young Scholar (2015A030306029).}
\\
\author{Chunhua Jin$^a$,  Yifu Wang$^b$, Jingxue Yin$^a$
\thanks{
Corresponding author. Email: {\tt yjx@scnu.edu.cn, jinchhua@126.com}}
\\
\small \it{$^a$School of Mathematical Sciences, South China Normal University, }
\\
\small \it{Guangzhou, 510631, China}
\\
\small \it{$^b$
School of Mathematics and Statistics, Beijing Institute of Technology, }
\\
\small \it{Beijing 100081, China}
}}

\date{}

\maketitle

\begin{abstract}
In this paper, we study a nutrient-taxis model with porous medium slow diffusion
\begin{align*}
\left\{
\begin{aligned}
&u_t=\Delta u^m-\chi\nabla\cdot(u\nabla v)+\xi uv-\rho u,
\\
&v_t-\Delta v=-vu+\mu v(1-v),
\end{aligned}\right.
\end{align*}
in a bounded domain $\Omega\subset \mathbb R^3$ with zero-flux boundary condition.
It is shown that for any $m>\frac{11}4-\sqrt 3$,
the problem admits a global weak solution for any large initial datum. We divide the study into three cases,
(i) $\xi\mu=0, \rho\ge 0$;   (ii) $\xi\mu\rho>0$;  (iii)  $\xi\mu>0$, $\rho=0$.
In particular, for Case (i) and Case (ii), the global solutions are uniformly bounded.
Subsequently, the large time behavior of these global bounded solutions are also discussed.
At last, we also extend the results to the coupled chemotaxis-Stokes system.
Important progresses for chemotaxis-Stokes system with $m>\frac 76$, $m>\frac 87$ and $m>\frac 98$ have been
carried out respectively by \cite{W2, TW2, W3}, but leave a gap for $1<m\le \frac98$. Our result
for chemotaxis-Stokes system
supplements part of the gap $(\frac{11}4-\sqrt 3, \frac 98)$. Here $\frac{11}4-\sqrt 3\approx 1.018$.
\end{abstract}


\section{Introduction}
In this paper, we consider the following nutrient-taxis model involving food-supported proliferation
\begin{align}
\label{1-1}\left\{
\begin{aligned}
&u_t=\Delta u^m-\chi\nabla\cdot(u\nabla v)+\xi uv-\rho u, \    (x,t) \in Q,
\\
&v_t-\Delta v=-vu+\mu v(1-v), \    (x,t) \in Q,
\\
&\left.(\nabla u^m-\chi u\nabla v)\cdot{\bf n}\right|_{\partial\Omega}
=\left.\frac{\partial v}{\partial {\bf n}}\right|_{\partial\Omega}=0,
\\
& u(x,0)=u_0(x),  v(x,0)=v_0(x), \quad  x\in\Omega,
\end{aligned}\right.
\end{align}
where $m>1$, $Q=\Omega\times \mathbb R^+$,
$\Omega\subset \mathbb R^3$ is a bounded domain, and the boundary $\partial\Omega$ is appropriately smooth,
$u$, $v$ represent the bacteria cell density,
the concentration of nutrient respectively, $\chi>0$
is the sensitivity coefficient of aggregation induced by the  concentration changes of nutrient,
the appearance of $\xi uv$ implies that the cell proliferation relies on the availability of nutrient resource $v$,
and $\xi\ge 0$ is the conversion rate (growth yield)
of consumed nutrient to bacterial growth, $-\rho u$ ($\rho\ge 0$) is the linear degradation of the  bacteria cells,
$-vu$ and $\mu v(1-v)$ with $\mu\ge 0$  represent the consumption and reproduction of nutrients, respectively.
In addition to the above biological explanation, this model is often used by some authors
to describe the prey-taxis phenomenon involving Lotka-Volterra type interaction,
see for example \cite{AB, LH, KO}.

Colonies of bacteria growing on the surface of thin agar plates show varieties of
morphological patterns in response to surrounding environmental conditions, such as the nutrient concentration,
the solidity of an agar medium and temperature. Based on experimental observations,
Kawasaki et al. \cite{KMM} proposed the following reaction-diffusion model
for bacterial aggregation patterns on the surface of thin agar plates
$$\left\{
\begin{aligned}
&u_t=\nabla\cdot(D_u\nabla u)+\xi f(u,v),
\\
&v_t=D_v\Delta v-f(u,v),
\end{aligned}\right.
$$
where $D_u$ and $D_v$ are the diffusion coefficients of the
bacterial cells and nutrient, respectively.  In recent years, bacterial chemotaxis
has attracted much attention due to its critical role in pattern formation.
To explore the aggregation patterns caused by such chemotactic  mechanism,
Leyva et. al \cite{LM} took nutrient chemotactic term into the above model,
and developed the following model
$$\left\{
\begin{aligned}
&u_t=\nabla\cdot(uv\nabla u)-\nabla\cdot(u^2v\nabla v)+uv,
\\
&v_t=\Delta v-uv.
\end{aligned}\right.
$$
The key feature of this model is the choice of the nonlinear diffusion coefficient $D_u$, which depends on both $u$ and $v$,
that is, the bacteria are immotile when
either  $u$ or $v$ are low and become active as $u$ or $v$ increase.
Recently, Winkler \cite{W} considered a simplified form of this model, that is
\begin{equation}\label{1-2}\left\{
\begin{aligned}
&u_t=\Delta u-\chi\nabla\cdot(u\nabla v)+\xi uv-\rho u,
\\
&v_t=\Delta v- uv+\mu v(1-\alpha v),
\end{aligned}\right.
\end{equation}
in which, the author included the possibility of linear degradation in the cell population,
and the reproduction of chemoattractant through either linear or logistic mechanism.
For this model, Winkler obtained the existence of global weak solutions, and further proved that
under some assumptions on these coefficients, each of these solutions becomes eventually
smooth and stabilizes toward a spatially homogeneous equilibrium.
When $\xi=\rho=\mu=0$, this model \eqref{1-2} is reduced to the following form
$$\left\{
\begin{aligned}
&u_t=\Delta u-\chi\nabla\cdot(u\nabla v),
\\
&v_t=\Delta v- uv.
\end{aligned}\right.
$$
In 2012, Tao and Winkler \cite{TW} showed the existence of  global weak solutions in three dimensions,
for which, they also proved that,
after some time $T$, these weak solutions become smooth and go to constant equilibria in the large time limit.
While if the cell mobility is described by a nonlinear function of the cells density, for example,
the porous medium diffusion, then the model becomes
\begin{equation*}\left\{
\begin{aligned}
&u_t=\Delta u^m-\chi\nabla\cdot(u\nabla v),
\\
&v_t=\Delta v- uv.
\end{aligned}\right.
\end{equation*}
If the fluid velocity is considered into this model, the system \eqref{1-3} becomes the classical
chemotaxis-Stokes system,
\begin{equation}\label{1-3}\left\{
\begin{aligned}
&n_t+u\cdot\nabla n=\Delta n^m-\chi\nabla\cdot(n\nabla c),
\\
&c_t++u\cdot\nabla c-\Delta c=-cn,
\\
&u_t+\nabla P=\Delta u+n\nabla\varphi,
\\
& {\rm div} u=0,
\end{aligned}\right.
\end{equation}
This model is introduced by Tuval, Goldstein, et.al \cite{TC} in 2005.
which describes the dynamics of bacterial swimming and oxygen transport near
contact lines. Since then, this model has been studied by many researchers.
For the two dimensional case of \eqref{1-3}, the global solvability and boundedness of weak solutions are established completely
for any $m>1$ in \cite{TW1}.
While in three dimensional space, the research of \eqref{1-3} is rather tortuous. The first effort to this 3-D problem is due to the
work by Di Francesco et al. \cite{DLM}, in which, they obtained the existence of global bounded
weak solutions for $m$ in some finite interval, namely  $m\in \left(\frac{7+\sqrt{217}}{12}, 2\right]$ (approximating to  $(1.8109, 2]$);
It was Tao and Winkler \cite{TW2}, in 2013, who established the  global existence of locally
bounded weak solutions with $m$ belonging to the infinite interval $(\frac87, +\infty)$.
Afterwards, Winkler \cite{W2} supplemented the uniform boundedness of solutions
for the case $m>\frac 76$; Recently, Winkler \cite{W3} further improved this result to the case $m>\frac 98$ with $\Omega$ being
a convex domain. However, as mentioned by Winkler \cite{W3}, the question of identifying an optimal condition on
$m\ge 1$ ensuring global boundedness in the three-dimensional version of \eqref{1-3} remains an open challenge.

In the present paper, we  first pay our attention to the global existence  and uniform boundedness  of weak solutions
for the system \eqref{1-1}.  We divide the research into three cases
according to the nonnegative coefficients $\rho, \xi, \mu$,  that is
$$
(i) \ \xi\mu=0, \rho\ge 0;   \qquad (ii) \  \xi\mu\rho>0; \qquad (iii)  \  \xi\mu>0, \rho=0.
$$
We show that for any $m>1$, this problem admits a global  weak solution
for any large initial datum and any nonnegative coefficients $\rho, \xi, \mu$. In particular,
the solution is uniformly bounded for Cases (i) and (ii), while for Case (iii), the solution
is just locally bounded on time $t$ since the $L^1$-norm of $u$ depends on $t$.

Throughout this paper, we assume that

\begin{align*}
\qquad\qquad\qquad\qquad
\left\{
\begin{aligned}
& u_0\in L^\infty(\Omega), \nabla u_0^m\in L^{2}(\Omega), v_0\in W^{2, \infty}(\Omega),
\\
& u_0, v_0\ge 0,
\\
&\partial\Omega\in C^{2,\alpha}.
\end{aligned}\right. \hspace {5cm} (H_1)
\end{align*}

In what follows, we give the existence results. For Cases (i) and (ii),  we have
\begin{theorem}
\label{thm-1}
Assume $(H_1)$,  $m>\frac{11}4-\sqrt 3$.  If (i) $\xi\mu=0$ with $\rho\ge 0$,  or (ii) $\xi\mu\rho>0$,
the problem \eqref{1-1} admits a nonnegative global bounded weak solution $(u,v)$
with  $u\in\mathcal X_1, v\in \mathcal X_2$,
where
$$
\mathcal X_1=\{u\in L^\infty(\Omega\times R^+); \nabla u^m\in L^\infty((0,\infty);L^2(\Omega)),
\left(u^{\frac{m+1}2}\right)_t, \nabla u^{\frac{m+1}2} \in L^2_{loc}([0,\infty); L^2(\Omega))\},
$$
$$
\mathcal X_2=\{v\in L^\infty((0,\infty); W^{1,\infty}(\Omega)); v_t, \Delta v\in  L^p_{loc}([0,\infty);L^p(\Omega))\ \  \text{for any }\  p>1\},
$$
such that
\begin{equation}
\label{1-4}
\sup_{t\in(0,+\infty)}\left(\|u(\cdot,t)\|_{L^\infty}+\|v(\cdot,t)\|_{W^{1,\infty}}\right)\le M_1,
\end{equation}
\begin{equation}
\label{1-5}
\sup_{t\in (0,+\infty)}\int_\Omega|\nabla u^m|^2 dx+
\sup_{t\in (0,+\infty)}\|u^{\frac{m+1}2}\|_{W_2^{1,1}(Q_1(t))}
\le M_2,
\end{equation}
\begin{equation}
\label{1-6}
\sup_{t\in (0,+\infty)}\|v\|_{W_p^{2,1}(Q_1(t))}
\le M_3 \quad \text{for any }\  p>1.
\end{equation}
Here $Q_1(t)=\Omega\times (t, t+1)$,
$M_i$ $(i=1,2,3)$ are constants depending only on $\xi$, $\chi$, $\rho$, $\mu$, $\Omega$, $u_0$, $v_0$.
\end{theorem}

For the case (iii), the global solution  of  \eqref{1-1} is
locally bounded, which can be stated as follows.
\begin{theorem}
\label{thm-2}
Assume $(H_1)$,  $m>\frac{11}4-\sqrt 3$.  If (iii) $\xi\mu>0$,  and $\rho=0$,
the problem \eqref{1-1} admits a nonnegative local bounded weak solution $(u,v)$
with  $u\in\tilde{\mathcal  X_1}, v\in  \tilde{\mathcal X_2}$,
where
$$
\tilde{\mathcal  X_1}=\{u\in L^\infty_{loc}([0,+\infty); L^\infty(\Omega) );\nabla u^m\in L_{loc}^\infty([0,+\infty);L^2(\Omega)),
\left(u^{\frac{m+1}2}\right)_t, \nabla u^{\frac{m+1}2} \in L_{loc}^2([0,+\infty); L^2(\Omega))\},
$$
$$
\tilde{\mathcal  X_2}=\{v\in L_{loc}^\infty([0,+\infty); W^{1,\infty}(\Omega)); v_t, \Delta v\in  L^p_{loc}([0,+\infty);L^p(\Omega))
\quad \text{for any }\  p>1\},
$$
such that for any $T>0$,
\begin{equation}
\label{1-7}
\sup_{t\in(0, T)}\left(\|u(\cdot,t)\|_{L^\infty}+\|v(\cdot,t)\|_{W^{1,\infty}}\right)\le \tilde M_1(T),
\end{equation}
\begin{equation}
\label{1-8}
\sup_{t\in (0, T)}\int_\Omega|\nabla u^m|^2 dx+
\|u^{\frac{m+1}2}\|_{W_2^{1,1}(Q_T)}
\le \tilde M_2(T),
\end{equation}
\begin{equation}
\label{1-9}
\|v\|_{W_p^{2,1}(Q_T)}
\le \tilde M_3(T) \quad \text{for any }\  p>1.
\end{equation}
Here $Q_T=\Omega\times (0, T)$, $\tilde M_i(T)$ $(i=1,2,3)$ are constants depending only on $T$,  $\xi$, $\chi$, $\rho$, $\mu$, $\Omega$, $u_0$, $v_0$.
\end{theorem}

On the basis of establishing the global solvability, we further consider the large time
behavior of the global solutions. We only consider the cases (i) and (ii), since the global solutions are bounded uniformly
for the two cases.
\begin{theorem}
\label{thm-3}
Assume $(H_1)$, $m>\frac{11}4-\sqrt 3$, and $u_0\not\equiv 0$, $v_0\not\equiv 0$.
Let $(u,v)$ be the global bounded solution obtained above. Then we have

{\rm (i)} When $\mu=0$, $\rho=0$, $\xi\ge 0$,  then
\begin{align}
\label{1-10}
\lim_{t\to\infty}\|v\|_{L^\infty}=0, \qquad \lim_{t\to\infty}\|u-A\|_{L^p}=0 \quad\text{for any $p>1$},
\end{align}
where $A=\frac{1}{|\Omega|}\int_\Omega(u_0+\xi v_0)dx>0$;

{\rm (ii)} When $\mu=0$, $\rho>0$, $\xi\ge 0$,  there exists a constant $B$ with $0<B<\frac{1}{|\Omega|}\int_\Omega v_0 dx$, such that
\begin{align}
\label{1-11}
\lim_{t\to\infty}\|v-B\|_{L^\infty}=0,  \qquad \lim_{t\to\infty}\|u\|_{L^\infty}=0;
\end{align}

{\rm (iii)}  When $\mu>0$, $\rho>0$, $0\le \xi< \rho$,  then
\begin{align}
\label{1-12}
\lim_{t\to\infty}\|v-1\|_{L^\infty}=0,  \qquad \lim_{t\to\infty}\|u\|_{L^\infty}=0.
\end{align}

\end{theorem}

At last, we also consider the following typical chemotaxis-Stokes system
\begin{align}
\label{1-13}\left\{
\begin{aligned}
&n_t+u\cdot\nabla n=\Delta n^m-\chi\nabla\cdot(n\nabla c), \    (x,t) \in Q,
\\
&c_t++u\cdot\nabla c-\Delta c=-cn, \   (x,t) \in Q,
\\
&u_t+\nabla P=\Delta u+n\nabla\varphi, \    (x,t) \in Q,
\\
& {\rm div} u=0, \    (x,t) \in Q,
\\
&\left.(\nabla n^m-\chi n\nabla c)\cdot{\nu}\right|_{\partial\Omega}
=\left.\frac{\partial c}{\partial {\nu}}\right|_{\partial\Omega}=0,   \quad u|_{\partial\Omega}=0,
\\
& n(x,0)=n_0(x),  c(x,0)=v_0(x), u(x,0)=u_0, \quad  x\in\Omega.
\end{aligned}\right.
\end{align}
Our research for this system supplements part of the gap of \cite{TW2, W2, W3} for $m\in(\frac{11}4-\sqrt 3,\frac98]$.

Assumptions:
\begin{align*}
\qquad\qquad\qquad\qquad
\left\{
\begin{aligned}
& n_0\in L^\infty(\Omega), \nabla n_0^m\in L^{2}(\Omega), c_0\in W^{2, \infty}(\Omega),
\\
&u_0\in L^\infty(\Omega), A^\beta u_0\in  L^2(\Omega), \text{for any $\beta\in(\frac34, 1)$},
{\rm div}u_0=0, \varphi\in W^{1,\infty}(\Omega),
\\
& n_0, c_0\ge 0,
\\
&\partial\Omega\in C^{2,\alpha}.
\end{aligned}\right.\quad (H_2)
\end{align*}
The main results:
\begin{theorem}
\label{thm-4}
Assume $(H_2)$,  $m>\frac{11}4-\sqrt 3$.
Then the problem \eqref{1-13} admits a  global bounded weak solution $(n,c,u,\pi)$
with  $n\in\mathcal X_1, c\in \mathcal X_2$, $u\in \mathcal X_3$, $\pi\in  \mathcal X_4$,
where
$$
\mathcal X_1=\{n\in L^\infty(\Omega\times \mathbb R^+); \nabla n^m\in L^\infty(\mathbb R^+;L^2(\Omega)),
\left(n^{\frac{m+1}2}\right)_t, \nabla n^{\frac{m+1}2} \in L^2_{loc}([0,\infty); L^2(\Omega))\},
$$
$$
\mathcal X_2=\{c\in L^\infty(\mathbb R^+; W^{1,\infty}(\Omega)); c_t, \Delta c\in  L^p_{loc}([0,\infty);L^p(\Omega))
\  \text{for any }\  p>1\},
$$
$$
\mathcal X_3=\left\{u\in  L^\infty(Q); \nabla u, \Delta u, u_t\in L^2_{loc}(\mathbb R^+;L^2(\Omega)),
A^\beta u\in L^\infty(\mathbb R^+;
L^2(\Omega)) \ \forall \ \beta\in(\frac34,1) \right\},
$$
$$
\mathcal X_4=\{\pi; \nabla\pi\in L^2_{loc}(\mathbb R^+;L^2(\Omega))\},
$$
such that
\begin{align}
\label{1-14}
\sup_{t\in(0,+\infty)}\left(\|n(\cdot,t)\|_{L^\infty}+\|c(\cdot,t)\|_{W^{1,\infty}}
+\|\nabla n^m(\cdot,t)\|_{L^2}+\|u(\cdot,t)\|_{L^\infty}+\|A^\beta u(\cdot,t)\|_{L^2}
\right)\le M_1,
\\
\label{1-15}
\sup_{t\in (0,+\infty)}\left(\|n^{\frac{m+1}2}\|_{W_2^{1,1}(Q_1(t))}+
\|n\|_{W_p^{2,1}(Q_1(t))}+\|u\|_{W_2^{2,1}(Q_1(t))}\right)\le M_2, \ \text{for any }\  p>1.
\end{align}
Here $Q_1(t)=\Omega\times (t, t+1)$,
$M_i$ $(i=1,2)$ are constants depending only on $\xi$, $\Omega$, $n_0$, $c_0$, $u_0$.
\end{theorem}
The stability of solutions for system \ref{1-13} had been established in \cite{W3},
that is the solutions go to the spatially homogeneous steady state $(\overline n_0, 0, 0)$
in the large time limit.

\section{Preliminaries}

We first give some notations, which will be used throughout this paper.

{\bf Notations:}
$\|\cdot\|_{L^p}=\|\cdot\|_{L^p(\Omega)}$,
$Q_{1}(t)=\Omega\times(t, t+1)$,  $Q_T:=Q_T(0)=\Omega\times (0, T)$.

Next, we give the definition of weak solutions.
\begin{definition}
$(u,v)$ is called a weak solution of \eqref{1-1}, if $u\ge 0, v\ge 0$ and
$u\in \mathcal X$, $v\in W_2^{1,0}(Q_T)$ for any $T>0$, such that
\begin{align*}
&-\iint_{Q_T}u\varphi_t dxdt-\int_\Omega u(x,0)\varphi(x,0)dx+\iint_{Q_T}(\nabla u^m-\chi
u\nabla v)\cdot\nabla\varphi dxdt
=\iint_{Q_T}(\xi uv-\rho u)\varphi dxdt,
\\
&-\iint_{Q_T}v\varphi_tdxdt-\int_\Omega v(x,0)\varphi(x,0)dx+\iint_{Q_T}\nabla v\cdot\nabla\varphi dxdt+\iint_{Q_T}
uv\varphi dxdt=\mu\iint_{Q_T} v(1-v)\varphi dxdt,
\end{align*}
for any $\varphi\in C^\infty(\overline Q_T)$ with  $\varphi(x,T)=0$,
where $\mathcal X=\{u\in L^2(Q_T); \nabla u^m\in L^2(Q_T)\}$.
\end{definition}

Before going further, we list some important lemmas, which will be used throughout this paper.
\begin{lemma}
\label{lem2-1}
Assume that $u\in L^p$, $v\in L^q$, $w\in L^r$, then
$$
\|uvw\|_{L^1}\le \|u\|_{L^p}\|v\|_{L^q}\|w\|_{L^r}\le \varepsilon_1\|u\|_{L^p}^p+\varepsilon_2\|u\|_{L^q}^q
+C(\varepsilon_1, \varepsilon_2)\|w\|_{L^r}^r.
$$
where $p,q,r>1$, $\frac1p+\frac1q+\frac 1r=1$, $\varepsilon_1$, $\varepsilon_2$ are two  arbitrarily small constants,
and $C(\varepsilon_1, \varepsilon_2)$ is a constant depends on $\varepsilon_1$, $\varepsilon_2$.
\end{lemma}
Next, by \cite{J1, J2}, we give the following two lemmas.
\begin{lemma}
\label{lem2-2}
Let $T>0$, $\tau\in (0, T)$, $\sigma\ge 0$, $a>0$, $b\ge 0$, and suppose that $f: [0, T)\to [0, \infty)$
is absolutely continuous, and satisfies
\begin{equation}
\label{2-1}
f'(t)+af^{1+\sigma}(t)\le h(t), t\in \mathbb R,
\end{equation}
where $h\ge 0$, $h(t)\in L^1_{loc}([0, T))$ and
$$
\int_{t-\tau}^t h(s)ds\le b, \ \text{for all}\ t\in[\tau, T).
$$
Then
\begin{equation}
\label{2-2}
\sup_{t\in(0, T)} f(t) +a\sup_{t\in(\tau, T)}\int_{t-\tau}^t f^{1+\sigma}(s) ds
\le b+2\max\{f(0)+b+a\tau, \frac{b}{a\tau}+1+2b+2a\tau\}.
\end{equation}
\end{lemma}
\begin{lemma}
\label{lem2-3}
Assume that $u_0\in W^{2,p}(\Omega)$, and  $f\in L_{loc}^p((0,+\infty);L^p(\Omega))$
with
$$
\sup_{t\in(\tau, +\infty)}\int_{t-\tau}^{t}\|f\|_{L^p}^pds\le A,
$$
where $\tau>0$ is a fixed constant.
Then the following problem
\begin{align}
\label{LP}\left\{
\begin{aligned}
&u_t-\Delta u+u=f(x,t),
\\
&\left.\frac{\partial u}{\partial{\bf n}}\right|_{\partial\Omega}=0,
\\
&u(x,0)=u_0(x)
\end{aligned}\right.
\end{align}
admits a unique solution $u$ with $u\in  L_{loc}^p((0,+\infty);W^{2,p}(\Omega))$, $u_t\in L_{loc}^p((0,+\infty);L^p(\Omega))$
with
\begin{equation}
\label{LPE}
\sup_{t\in(\tau, +\infty)}\int_{t-\tau}^{t}(\|u\|_{W^{2,p}}^p+\|u_t\|_{L^p}^p)ds\le
AM\frac{e^{p\tau}}{e^{\frac p2\tau}-1}+Me^{\frac p2\tau}\|u_0\|_{W^{2,p}}^p,
\end{equation}
where $M$ is a constant independent of $\tau$.
\end{lemma}

By \cite{W1}, We also have the following lemma.
\begin{lemma}
\label{lem2-4}
Suppose that $h\in C^2(\mathbb R)$, then for all $\varphi\in C^2(\overline \Omega)$
fulfilling $\frac{\partial\varphi}{\partial {\bf n}}=0$ on $\partial\Omega$, we have
\begin{align}
&\int_\Omega h'(\varphi)|\nabla\varphi|^2\Delta\varphi dx+\frac23\int_\Omega
h(\varphi)|\Delta\varphi|^2 dx\nonumber
\\
\label{2-5}
=&\frac23\int_\Omega h(\varphi)|D^2\varphi|^2 dx
-\frac13\int_\Omega h''(\varphi)|\nabla\varphi|^4dx-\frac13\int_{\partial\Omega}
h(\varphi)\frac{\partial|\nabla\varphi|^2}{\partial {\bf n}}ds,
\end{align}
and
\begin{equation}
\label{2-6}
\int_\Omega\frac{|\nabla \varphi|^4}{\varphi^3}dx\le (2+\sqrt N)^2\int_\Omega \varphi|D^2\ln \varphi|^2dx.
\end{equation}
\end{lemma}

By \cite{MS}, we have
\begin{lemma}
\label{lem2-5}
Assume that $\Omega$ is bounded and let $\omega\in C^2(\overline\Omega)$ satisfy
$\frac{\partial\omega}{\partial\nu}\Big|_{\partial\Omega}=0$. Then we have
$$
\frac{\partial|\nabla\omega|^2}{\partial\nu}\le 2\kappa|\nabla\omega|^2 \quad \text{on} \ \partial\Omega,
$$
where  $\kappa>0$ is an upper bound for the curvatures of $\Omega$.
\end{lemma}

\section{Boundedness and Global Existence of Weak Solutions}
We first consider the approximate problems given by
\begin{align}
\label{3-1}\left\{
\begin{aligned}
&u_{\varepsilon t}=\Delta (\varepsilon u_{\varepsilon} +u_{\varepsilon}^m)-\chi
\nabla\cdot\left(u_{\varepsilon}\nabla v_{\varepsilon}\right)+\xi u_{\varepsilon}v_{\varepsilon}
-\rho u_{\varepsilon}-\varepsilon u_{\varepsilon}^2,  \  (x,t) \in Q,
\\
&v_{\varepsilon t}-\Delta v_{\varepsilon}=-v_{\varepsilon}u_{\varepsilon}
+\mu v_{\varepsilon}(1-v_{\varepsilon}), \    (x,t) \in Q,
\\
&\left.\frac{\partial u_{\varepsilon}}{\partial {\bf n}}\right|_{\partial\Omega}
=\left.\frac{\partial v_{\varepsilon}}{\partial {\bf n}}\right|_{\partial\Omega}=0,
\\
& u_{\varepsilon}(x,0)=u_{\varepsilon 0}(x),  v_{\varepsilon}(x,0)=v_{\varepsilon 0}(x), \quad  x\in\Omega.
\end{aligned}\right.
\end{align}
where $u_{\varepsilon 0}, v_{\varepsilon 0} \in C^{2+\alpha}(\overline\Omega)$ with
$\left.\frac{\partial u_{\varepsilon 0}}{\partial{\bf n}}\right|_{\partial\Omega}=
0, \left.\frac{\partial v_{\varepsilon 0}}{\partial{\bf n}}\right|_{\partial\Omega}=0$,
$\|u_{\varepsilon 0}\|_{L^\infty}+\|\nabla u_{{\varepsilon 0}}^m\|_{L^2}+\|v_{\varepsilon 0}\|_{W^{2,\infty}}
\le 2(\|u_0\|_{L^\infty}+\|\nabla u_{0}^m\|_{L^2}+\|v_0\|_{W^{2,\infty}})=M_0$, and
$$
u_{\varepsilon 0}\to u_0, v_{\varepsilon 0}\to v_0,\quad \text{uniformly}.
$$
According to the arguments in \cite{J2},  each of these problems  is globally
solvable in the classical sense.

\begin{lemma}
\label{lem3-1}
Assume that $m>1$, then for any $\varepsilon>0$, the problem \eqref{3-1} admits a unique nonnegative classical solution
$(u_\varepsilon, v_\varepsilon)\in C^{2+\alpha, 1+\frac\alpha 2}(\overline\Omega\times[0,+\infty))$.
\end{lemma}

Using this lemma, we show the global existence of solutions for the problem \eqref{1-1}.
For this purpose, we show some a prior estimates of solutions.
In what follows, we let $C$, $C_i$, $\tilde C$ denote some different constants,
which are independent of $\varepsilon$, and if no special explanation, these constants
depend at most on $\Omega$, $\chi$, $\xi$, $\rho$, $\mu$,  $u_0$, $v_0$.

We first give the following lemma
\begin{lemma}
\label{lem3-2}
Assume that $m>1$.
Let $(u_\varepsilon, v_\varepsilon)$ be the solution of \eqref{3-1}, then we have
\begin{align}
\label{3-2}
&\sup_{t\in(0,\infty)}\|v_\varepsilon(\cdot,t)\|_{L^\infty}\le C_1, &&
\\
\label{3-3}
&\sup_{t\in(0,\infty)}\int_\Omega u_\varepsilon dx\le C_2, && \text{if $\xi\mu=0$, or
$\xi\mu\rho>0$},
\\
\label{3-4}
&\int_\Omega u_\varepsilon(x,t) dx\le C_3(1+t), && \text{if $\xi\mu>0$ and $\rho=0$},
\end{align}
where $C_1$, $C_2$ and $C_3$ are independent of $\varepsilon$.
\end{lemma}

{\it\bfseries Proof.} Firstly, consider the initial value problem of the following ODE
$$
\left\{
\begin{aligned}
&y'(t)=\mu y(1-y),  \quad t>0,
\\
&y(0)=M_0.
\end{aligned}\right.
$$
It is easy to obtain that $0\le y\le \max\{1, M_0\}$.
By comparison lemma, we obtain that for any $t>0$,
$$
\|v(\cdot,t)\|_{L^\infty}\le y(t)\le\max\{1,\|v_0\|_{L^\infty}\}.
$$
Integrating the first equation and the second equation respectively,
we obtain
$$
\frac{d}{dt}\int_\Omega u_{\varepsilon}dx+\rho\int_\Omega u_{\varepsilon}dx
+\varepsilon\int_\Omega u_{\varepsilon}^2dx
=\xi\int_\Omega u_{\varepsilon}v_{\varepsilon}dx
$$
and
$$
\xi\frac{d}{dt}\int_\Omega v_{\varepsilon}dx=-\xi\int_\Omega v_{\varepsilon}u_{\varepsilon}dx
+\xi\mu\int_\Omega v_{\varepsilon}(1-v_{\varepsilon})dx.
$$
Adding up the above two equalities, we obtain
$$
\frac{d}{dt}\int_\Omega (u_{\varepsilon}+\xi v_{\varepsilon})dx
+\rho\int_\Omega (u_{\varepsilon}+\xi v_{\varepsilon})dx
\le \rho\xi\int_\Omega v_{\varepsilon}dx+\xi\mu\int_\Omega v_{\varepsilon}(1-v_{\varepsilon})dx.
$$
When $\xi\mu=0$, it is easy to obtain
$$
\int_\Omega (u_{\varepsilon}+\xi v_{\varepsilon})dx\le \int_\Omega (u_{\varepsilon 0}+\xi v_{\varepsilon 0})dx.
$$
When $\xi\mu>0$,  we obtain
$$
\frac{d}{dt}\int_\Omega (u_{\varepsilon}+\xi v_{\varepsilon})dx
+\rho\int_\Omega (u_{\varepsilon}+\xi v_{\varepsilon})dx
\le C,
$$
if $\rho>0$, by a direct calculation, we obtain \eqref{3-3};
while if $\rho=0$, then
$$
\int_\Omega (u_{\varepsilon}+\xi v_{\varepsilon})dx\le \int_\Omega (u_{\varepsilon 0}+\xi v_{\varepsilon 0})dx+Ct.
$$
The proof is complete. \hfill $\Box$

From this lemma, we see that for the cases (i) $\xi\mu=0$, $\rho\ge 0$  and  (ii)  $\xi\mu\rho>0$,
the $L^1$-norm of $u_\varepsilon(\cdot, t)$  is uniformly bounded on $t$.
However, for the case (iii) $\xi\mu>0$ and $\rho=0$, the $L^1$-norm of $u_\varepsilon(\cdot, t)$
depends on $t$.  In what follows, we only show the energy estimates independent of time $t$ for the cases (i) and (ii).
For the case (iii), the similar energy estimates also hold,  but depend on time $t$.

\begin{lemma}
\label{lem3-3}
Assume that $m>1$, let $(u_\varepsilon, v_\varepsilon)$ be the solution of \eqref{3-1}.
Then for Cases (i) and (ii), we have
\begin{align}
&\sup_{t\in(0,+\infty)}\int_\Omega\left(\frac{|\nabla v_\varepsilon|^2}{v_\varepsilon}+
u_{\varepsilon}\ln u_{\varepsilon}\right)dx+\sup_{t\in(0,+\infty)}\int_{t}^{t+1}
\int_\Omega\left(\varepsilon u_{\varepsilon}^2\ln(1+ u_{\varepsilon})+
\varepsilon\frac{|\nabla u_{\varepsilon}|^2}{u_{\varepsilon}}
\right) dx
\nonumber
\\
\label{3-5}
&+\sup_{t\in(0,+\infty)}\int_{t}^{t+1}
\int_\Omega \left(v_{\varepsilon}|D^2\ln v_{\varepsilon}|^2 dx
+|\nabla u_{\varepsilon}^{\frac m2}|^2+\frac{u_\varepsilon}{v_\varepsilon}|\nabla v_\varepsilon|^2
+\frac{|\nabla v_{\varepsilon}|^4}{v_{\varepsilon}^3}+u_{\varepsilon}^{m+\frac 23}
\right)dx\le C,
\end{align}
where $C$ is independent of $\varepsilon$.
\end{lemma}

{\it\bfseries Proof.} Using the second equation of \eqref{3-1}, we obtain
\begin{align}
&\frac12\frac{d}{dt}\int_\Omega\frac{|\nabla v_\varepsilon|^2}{v_\varepsilon}dx
=\int_\Omega\frac{\nabla v_\varepsilon}{v_\varepsilon}\nabla v_{\varepsilon t}dx
-\frac12\int_\Omega\frac{|\nabla v_\varepsilon|^2}{v_\varepsilon^2}v_{\varepsilon t}dx
\nonumber
\\
&=-\int_\Omega v_{\varepsilon t}\left(\frac{\Delta v_{\varepsilon}}{v_{\varepsilon}}
-\frac{|\nabla v_\varepsilon|^2}{v_\varepsilon^2}\right)
-\frac12\int_\Omega\frac{|\nabla v_\varepsilon|^2}{v_\varepsilon^2}v_{\varepsilon t}dx
\nonumber
\\
&=-\int_\Omega \frac{\Delta v_{\varepsilon}}{v_{\varepsilon}} \left(
\Delta v_{\varepsilon}- v_{\varepsilon}u_{\varepsilon}
+\mu v_{\varepsilon}(1-v_{\varepsilon})\right)dx
+\frac12\int_\Omega \frac{|\nabla v_\varepsilon|^2}{v_\varepsilon^2}\left(
\Delta v_{\varepsilon}-v_{\varepsilon}u_{\varepsilon}
+\mu v_{\varepsilon}(1-v_{\varepsilon})\right)dx
\nonumber
\\
\label{3-6}
&=-\int_\Omega \frac{|\Delta v_{\varepsilon}|^2}{v_{\varepsilon}}
+\frac12\int_\Omega \frac{|\nabla v_\varepsilon|^2}{v_\varepsilon^2}\Delta v_{\varepsilon} dx
-\int_\Omega \left(\nabla u_{\varepsilon}\nabla v_{\varepsilon}
+\frac{3\mu}2|\nabla v_\varepsilon|^2+\frac{1}2\frac{u_\varepsilon}{v_\varepsilon}
|\nabla v_\varepsilon|^2-\frac{\mu}2\frac 1{v_\varepsilon}|\nabla v_\varepsilon|^2 dx\right).
\end{align}
Noticing that
\begin{align*}
\int_\Omega v_{\varepsilon}|D^2\ln v_{\varepsilon}|^2 dx &=\int_\Omega
\left(\frac{|D^2v_{\varepsilon}|^2}{v_{\varepsilon}}+\frac{|\nabla v_\varepsilon|^4}
{v_\varepsilon^3}-2\frac{\nabla v_\varepsilon \nabla^2
v_\varepsilon\nabla v_\varepsilon}{v_\varepsilon^2}\right) dx
\\
&=\int_\Omega
\left(\frac{|D^2v_{\varepsilon}|^2}{v_{\varepsilon}}+\frac{|\nabla v_\varepsilon|^4}
{v_\varepsilon^3}-\frac{\nabla v_\varepsilon\nabla |\nabla v_\varepsilon|^2}{v_\varepsilon^2}
\right)dx
\\
&=\int_\Omega
\left(\frac{|D^2v_{\varepsilon}|^2}{v_{\varepsilon}}+\frac{|\nabla v_\varepsilon|^4}
{v_\varepsilon^3}+\frac{|\nabla v_\varepsilon|^2}{v_\varepsilon^2}\Delta v_\varepsilon
-2\frac{|\nabla v_\varepsilon|^4}{v_\varepsilon^3}\right)dx
\\
&=\int_\Omega
\left(\frac{|D^2v_{\varepsilon}|^2}{v_{\varepsilon}}+\frac{|\nabla v_\varepsilon|^2}{v_\varepsilon^2}\Delta v_\varepsilon
-\frac{|\nabla v_\varepsilon|^4}{v_\varepsilon^3}\right)dx.
\end{align*}
By strong maximum principle, $v_\varepsilon(x,t)>0$ for any $x\in\Omega$, $t>0$.
Using Lemma \ref{lem2-4} and combining with the above equality,  we see that
\begin{align*}
&-\int_\Omega \frac{|\Delta v_{\varepsilon}|^2}{v_{\varepsilon}}
+\frac12\int_\Omega \frac{|\nabla v_\varepsilon|^2}{v_\varepsilon^2}\Delta v_{\varepsilon} dx
\nonumber
\\
=&-\int_\Omega\frac{|\nabla v_\varepsilon|^2}{v_\varepsilon^2}\Delta v_{\varepsilon} dx
-\int_\Omega\frac1{v_\varepsilon}|D^2v_\varepsilon|^2 +\frac{1}{v_\varepsilon^3}
|\nabla v_\varepsilon|^4 dx+\frac12\int_{\partial\Omega}\frac{1}{v_\varepsilon}
\frac{\partial}{\partial{\bf n}}|\nabla v_\varepsilon|^2 ds \nonumber
\\
=&-\int_\Omega v_{\varepsilon}|D^2\ln v_{\varepsilon}|^2 dx
+\frac12\int_{\partial\Omega}\frac{1}{v_\varepsilon}
\frac{\partial}{\partial{\bf n}}|\nabla v_\varepsilon|^2 ds.
\end{align*}
Substituting this inequality into \eqref{3-6}, we see that
\begin{align}
&\frac12\frac{d}{dt}\int_\Omega\frac{|\nabla v_\varepsilon|^2}{v_\varepsilon}dx
+\int_\Omega v_{\varepsilon}|D^2\ln v_{\varepsilon}|^2 dx+
\int_\Omega \left(\frac{3\mu}2|\nabla v_\varepsilon|^2
+\frac{1}2\frac{u_\varepsilon}{v_\varepsilon}|\nabla v_\varepsilon|^2\right)dx\nonumber
\\
\label{3-7}
&=\frac12\int_{\partial\Omega}\frac{1}{v_\varepsilon}
\frac{\partial}{\partial{\bf n}}|\nabla v_\varepsilon|^2 ds- \int_\Omega
\nabla u_{\varepsilon}\nabla v_{\varepsilon}dx
+\frac{\mu}2\int_\Omega\frac 1{v_\varepsilon}|\nabla v_\varepsilon|^2 dx.
\end{align}
Multiplying the first equation of \eqref{3-1} by $1+\ln u_{\varepsilon}$,  and integrating
it over $\Omega$,  we obtain
\begin{align}
&\frac{d}{dt}\int_\Omega u_{\varepsilon}\ln u_{\varepsilon}dx+\varepsilon\int_\Omega
\frac{|\nabla u_{\varepsilon}|^2}{u_{\varepsilon}}dx+
\frac{4}{m}\int_\Omega
|\nabla u_{\varepsilon}^{\frac m2}|^2dx+\int_\Omega\rho u_{\varepsilon}(1+\ln u_{\varepsilon})dx
+\varepsilon\int_\Omega u_{\varepsilon}^2(1+\ln u_{\varepsilon})dx\nonumber
\\
\label{3-8}
=&\chi\int_\Omega \nabla u_{\varepsilon}\nabla v_{\varepsilon} dx
+\xi\int_\Omega u_{\varepsilon} v_{\varepsilon} (1+\ln u_{\varepsilon}) dx,
\end{align}
using \eqref{3-2},  and it implies that
\begin{align}
&\frac{d}{dt}\int_\Omega u_{\varepsilon}\ln u_{\varepsilon}dx+\varepsilon\int_\Omega
\frac{|\nabla u_{\varepsilon}|^2}{u_{\varepsilon}}dx
+\frac{4}{m}\int_\Omega
|\nabla u_{\varepsilon}^{\frac m2}|^2dx+\int_\Omega\rho u_{\varepsilon}\ln (1+u_{\varepsilon})dx
+\varepsilon\int_\Omega u_{\varepsilon}^2\ln(1+ u_{\varepsilon})dx\nonumber
\\
\label{3-9}
\le &\chi\int_\Omega \nabla u_{\varepsilon}\nabla v_{\varepsilon} dx
+\eta\int_\Omega u_{\varepsilon}^{m+\frac23} dx +C_\eta
\end{align}
for any small $\eta>0$.
By Gagliardo-Nirenberg interpolation inequality,  we have
\begin{align}
\|u_{\varepsilon}\|_{L^{m+\frac23}}^{m+\frac23}=\|u_{\varepsilon}^{\frac m2}\|_{L^{2+\frac{4}{3m}}}^{2+\frac{4}{3m}}
&\le C_1\|u_{\varepsilon}^{\frac m2}\|_{L^{\frac 2m}}^{\frac{4}{3m}}\|\nabla u_{\varepsilon}^{\frac m2}\|_{L^2}^2
+C_2\|u_{\varepsilon}\|_{L^1}^{m+\frac23}\nonumber
\\
\label{3-10}
&=C_1\|u_{\varepsilon}\|_{L^1}^{\frac23}\|\nabla u_{\varepsilon}^{\frac m2}\|_{L^2}^2
+C_2\|u_{\varepsilon}\|_{L^1}^{m+\frac23}.
\end{align}
Substituting the above inequality into \eqref{3-9},  and using Lemma \ref{lem3-2}, we see that
\begin{align}
&\frac{d}{dt}\int_\Omega u_{\varepsilon}\ln u_{\varepsilon}dx+\varepsilon\int_\Omega
\frac{|\nabla u_{\varepsilon}|^2}{u_{\varepsilon}}dx
+\frac{2}{m}\int_\Omega
|\nabla u_{\varepsilon}^{\frac m2}|^2dx+\int_\Omega\rho u_{\varepsilon}\ln (1+u_{\varepsilon})dx
+\varepsilon\int_\Omega u_{\varepsilon}^2\ln(1+ u_{\varepsilon})dx\nonumber
\\
\label{3-11}
\le &\chi\int_\Omega \nabla u_{\varepsilon}\nabla v_{\varepsilon} dx+C.
\end{align}
Combining \eqref{3-11} with \eqref{3-7},  we obtain
\begin{align}
&\frac{d}{dt}\int_\Omega\left(\frac12\frac{|\nabla v_\varepsilon|^2}{v_\varepsilon}+
\frac{1}{\chi}u_{\varepsilon}\ln u_{\varepsilon}\right)dx
+\int_\Omega v_{\varepsilon}|D^2\ln v_{\varepsilon}|^2 dx+\frac{\varepsilon}{\chi}\int_\Omega
\frac{|\nabla u_{\varepsilon}|^2}{u_{\varepsilon}}dx+\frac{2}{m\chi}\int_\Omega
|\nabla u_{\varepsilon}^{\frac m2}|^2dx\nonumber
\\
&+\frac{1}{\chi}\int_\Omega\rho u_{\varepsilon}\ln (1+u_{\varepsilon})dx
+\frac{\varepsilon}{\chi}\int_\Omega u_{\varepsilon}^2\ln(1+ u_{\varepsilon})dx
+\int_\Omega \left(\frac{3\mu}2|\nabla v_\varepsilon|^2
+\frac{1}2\frac{u_\varepsilon}{v_\varepsilon}|\nabla v_\varepsilon|^2\right)dx\nonumber
\\
\label{3-12}
&\le\frac12\int_{\partial\Omega}\frac{1}{v_\varepsilon}
\frac{\partial}{\partial{\bf n}}|\nabla v_\varepsilon|^2 ds
+\frac{\mu}2\int_\Omega\frac 1{v_\varepsilon}|\nabla v_\varepsilon|^2 dx+C.
\end{align}
By the boundary trace embedding theorem \cite{A}  and Lemma \ref{lem2-5},  we see that
\begin{align}
&\frac12\int_{\partial\Omega}\frac{1}{v_\varepsilon}
\frac{\partial}{\partial{\bf n}}|\nabla v_\varepsilon|^2 ds\le \kappa
\int_{\partial\Omega}\frac{1}{v_\varepsilon}|\nabla v_\varepsilon|^2 ds=
\kappa\int_{\partial\Omega}|v_\varepsilon^{\frac12}\nabla\ln v_\varepsilon|^2 ds\nonumber
\\
\le &\eta_1\kappa\int_{\Omega}|D (v_\varepsilon^{\frac12}\nabla\ln v_\varepsilon)|^2 dx+
C_{\eta_1}\int_{\Omega}v_\varepsilon|\nabla\ln v_\varepsilon|^2 dx \nonumber
\\
\le &\eta_1\kappa\int_{\Omega}\left(\frac12 v_\varepsilon^{-\frac12}\nabla v_\varepsilon\nabla\ln v_\varepsilon
+v_\varepsilon^{\frac12}D^2\ln v_\varepsilon\right)^2 dx+
C_{\eta_1}\int_{\Omega}\frac{|\nabla v_\varepsilon|^2}{v_\varepsilon} dx\nonumber
\\
\le & 2\eta_1\kappa\int_{\Omega}
\left(v_\varepsilon |D^2\ln v_\varepsilon|^2+\frac14\frac{|\nabla v_\varepsilon|^4}{v_\varepsilon^3}\right) dx
+C_{\eta_1}\int_{\Omega}\frac{|\nabla v_\varepsilon|^2}{v_\varepsilon} dx\nonumber
\\
\label{3-13}
\le & 10\eta_1\kappa\int_{\Omega} v_\varepsilon |D^2\ln v_\varepsilon|^2 dx
+C_{\eta_1}\int_{\Omega}\frac{|\nabla v_\varepsilon|^2}{v_\varepsilon} dx
\end{align}
for any sufficiently small $\eta_1>0$.
By \eqref{3-2} and \eqref{2-6},  for any sufficiently small $\eta_2>0$,  we obtain
\begin{align}
\int_{\Omega}\frac{|\nabla v_\varepsilon|^2}{v_\varepsilon} dx &\le
\eta_2\int_{\Omega}\frac{|\nabla v_\varepsilon|^4}{v_\varepsilon^3} dx+C_{\eta_2}\nonumber
\\
\label{3-14}
&\le 16\eta_2\int_{\Omega} v_\varepsilon |D^2\ln v_\varepsilon|^2 dx+C_{\eta_2}.
\end{align}
Substituting \eqref{3-13}, \eqref{3-14}  into \eqref{3-12} gives
\begin{align}
&\frac{d}{dt}\int_\Omega\left(\frac12\frac{|\nabla v_\varepsilon|^2}{v_\varepsilon}+
\frac{1}{\chi}u_{\varepsilon}\ln u_{\varepsilon}\right)dx
+\frac12\int_\Omega v_{\varepsilon}|D^2\ln v_{\varepsilon}|^2 dx+\frac{\varepsilon}{\chi}\int_\Omega
\frac{|\nabla u_{\varepsilon}|^2}{u_{\varepsilon}}dx+\frac{2}{m\chi}\int_\Omega
|\nabla u_{\varepsilon}^{\frac m2}|^2dx\nonumber
\\
\label{3-15}
&+\frac{1}{\chi}\int_\Omega\rho u_{\varepsilon}\ln (1+u_{\varepsilon})dx
+\frac{\varepsilon}{\chi}\int_\Omega u_{\varepsilon}^2\ln(1+ u_{\varepsilon})dx
+\int_\Omega \left(\frac{3\mu}2|\nabla v_\varepsilon|^2
+\frac{1}2\frac{u_\varepsilon}{v_\varepsilon}|\nabla v_\varepsilon|^2\right)dx\le C.
\end{align}
By \eqref{2-6},   Lemma \ref{lem2-2}  and \eqref{3-10},  we  complete the proof. \hfill $\Box$

For the convenience of the proof. In what follows, we first separate out the case $m>2$.
By a direct calculation, it is easy to obtain the following lemma.
\begin{lemma}
\label{lem3-4}
Assume  (i) or (ii), and
let $(u_\varepsilon, v_\varepsilon)$ be the solution of \eqref{3-1}.  When $m>2$, for any $r>0$, we have
\begin{equation}
\label{3-16}
\sup_{t\in(0, +\infty)}\|u_\varepsilon(\cdot,t)\|_{L^{r+1}}^{r+1} +
\sup_{t\in(0, +\infty)}\int_t^{t+1}\|\nabla u_\varepsilon^{\frac{m+r}2}(\cdot,s)\|_{L^{2}}^{2} ds
\le C_r,
\end{equation}
and
\begin{equation}
\label{3-17}
\sup_{t\in(0, +\infty)}\|v_\varepsilon\|_{W_r^{2,1}Q_1(t)}\le \tilde C_r,
\end{equation}
where $C_r$, $\tilde C_r$ depend on $r$, and are independent of $\varepsilon$.
\end{lemma}

{\it\bfseries Proof. }  Multiplying the first equation of \eqref{3-1}  by $u_\varepsilon^r$ for any $r>0$,
then integrating it over $\Omega$, and using \eqref{3-2},  we obtain
\begin{align}
&\frac1{r+1}\frac{d}{dt}\int_\Omega u_\varepsilon^{r+1} dx+\varepsilon r
\int_\Omega u_\varepsilon^{r-1}|\nabla u_\varepsilon|^2 dx
+rm\int_\Omega u_\varepsilon^{m+r-2}|\nabla u_\varepsilon|^2 dx+\int_\Omega u_\varepsilon^{r+1} dx
+\varepsilon\int_\Omega u_\varepsilon^{r+2} dx\nonumber
\\
= & \xi\int_\Omega u_\varepsilon^{r+1}v_\varepsilon dx+r\chi\int_\Omega u_\varepsilon^r\nabla v_\varepsilon
\nabla u_\varepsilon dx+\int_\Omega u_\varepsilon^{r+1} dx\nonumber
\\
\label{3-18}
\le &C\int_\Omega u_\varepsilon^{r+1} dx+
\frac{rm}4\int_\Omega u_\varepsilon^{m+r-2}|\nabla u_\varepsilon|^2 dx+
C\int_\Omega u_\varepsilon^{r+2-m}|\nabla v_\varepsilon|^{2}dx.
\end{align}
By Gagliardo-Nirenberg interpolation inequality,  we see that
\begin{align*}
C\|u_\varepsilon\|_{L^{r+1}}^{r+1}&=C\|u_\varepsilon^{\frac{m+r}2}\|_{L^{\frac{2(r+1)}{m+r}}}^{\frac{2(r+1)}{m+r}}
\le C_1\|u_\varepsilon^{\frac{m+r}{2}}\|_{L^{\frac{2}{m+r}}}^{\frac{6m+4r-2}{(m+r)(3m+3r-1)}}
\|\nabla u_\varepsilon^{\frac{m+r}{2}}\|_{L^2}^{\frac{6r}{3m+3r-1}}+C_2\|u_\varepsilon\|_{L^1}^{r+1}
\\
&\le C_3(1+\|\nabla u_\varepsilon^{\frac{m+r}{2}}\|_{L^2}^{\frac{6r}{3m+3r-1}})
\le \frac{mr}{(m+r)^2}\|\nabla u_\varepsilon^{\frac{m+r}{2}}\|_{L^2}^2+C_4.
\end{align*}
Substituting the above inequality into \eqref{3-18}, we obtain
\begin{align}
&\frac1{r+1}\frac{d}{dt}\int_\Omega u_\varepsilon^{r+1} dx
+\frac{rm}2\int_\Omega u_\varepsilon^{m+r-2}|\nabla u_\varepsilon|^2 dx+\int_\Omega u_\varepsilon^{r+1} dx
\nonumber
\\
\label{3-19}
\le & C\int_\Omega u_\varepsilon^{r+2-m}|\nabla v_\varepsilon|^{2}dx+C_4.
\end{align}
By the above proof, it is easy to see that \eqref{3-19} hold for any $m>1$.  In what follows,
this inequality will be used often even for the case $m\le 2$,
and we don't prove it again.

Noticing  that when $m>2$, $r+3-m<r+1$, then by Gagliardo-Nirenberg interpolation inequality, we have
\begin{align*}
C\int_\Omega u_\varepsilon^{r+2-m}|\nabla v_\varepsilon|^{2}dx\le &C_1\int_\Omega u_\varepsilon^{r+3-m}dx+
C_2\int_\Omega |\nabla v_\varepsilon|^{2(r+3-m)}dx
\\
&\le \frac14 \int_\Omega u_\varepsilon^{r+1}dx+C_3
\|v_\varepsilon\|_{L^\infty}^{r+3-m}\int_\Omega |\Delta v_\varepsilon|^{r+3-m}dx+C_4.
\end{align*}
Substituting this inequality into \eqref{3-19}, we obtain
\begin{align}
&\frac1{r+1}\frac{d}{dt}\int_\Omega u_\varepsilon^{r+1} dx
+\frac{rm}2\int_\Omega u_\varepsilon^{m+r-2}|\nabla u_\varepsilon|^2 dx+\frac12\int_\Omega u_\varepsilon^{r+1} dx
\nonumber
\\
\label{3-20}
\le & C\int_\Omega|\Delta v_\varepsilon|^{r+3-m}dx+C_4.
\end{align}
As the application of Lemma \ref{lem2-3}, we have
\begin{align*}
&\sup_{t\in(0, +\infty)}\|u_\varepsilon(\cdot,t)\|_{L^{r+1}}^{r+1}+
\sup_{t\in(0, +\infty)}\int_t^{t+1}(\|\nabla u_\varepsilon^{\frac{m+r}2}(\cdot,s)\|_{L^{2}}^{2}+
\|u_\varepsilon(\cdot,s)\|_{L^{r+1}}^{r+1}) ds
\\
&\le C\sup_{t\in(0, +\infty)}\int_t^{t+1}\int_\Omega|\Delta v_\varepsilon|^{r+3-m}dxds+C
\\
&\le \tilde C\sup_{t\in(0, +\infty)}\int_t^{t+1}\int_\Omega u_\varepsilon^{r+3-m}dxds+\tilde C
\\
&\le \frac12 \sup_{t\in(0, +\infty)}\int_t^{t+1}\int_\Omega u_\varepsilon^{r+1}dxds+\hat C.
\end{align*}
Thus, \eqref{3-16} is proved, and \eqref{3-17} is a direct consequence of \eqref{3-16} and Lemma \ref{lem2-3}.
 \hfill $\Box$

In what follows, we construct  a sequence of iterations.
Taking $r=m-1$ in \eqref{3-19}, and noticing
that $u_\varepsilon|\nabla v_\varepsilon|^{2}\le \|v_\varepsilon\|_{L^\infty}
\frac{u_\varepsilon}{v_\varepsilon}|\nabla v_\varepsilon|^{2}$,
then combining with  \eqref{3-5} yields the following lemma, from which, we give the
initial  value of the sequence of iterations.

\begin{lemma}
\label{lem3-5}
Assume (i) or (ii), and $m>1$. Let $(u_\varepsilon, v_\varepsilon)$ be the solution of \eqref{3-1}. Then we have
\begin{equation}
\label{3-21}
\sup_{t\in(0,+\infty)}\|u_\varepsilon\|_{L^{m}}+\sup_{t\in(0,+\infty)}\int_t^{t+1}
\int_\Omega u_\varepsilon^{2m-3}|\nabla u_\varepsilon|^2 dxds\le C,
\end{equation}
where $C$ is independent of $\varepsilon$.
\end{lemma}

In order to make the iteration process proceed smoothly, we also give the following lemma

\begin{lemma}
\label{lem3-6}
Assume (i) or (ii), and $1<m\le 2$.
Let $(u_\varepsilon, v_\varepsilon)$ be the solution of \eqref{3-1}. If
$$
\sup_{t\in(0,+\infty)}\int_t^{t+1}
\int_\Omega u_\varepsilon^{\alpha}|\nabla u_\varepsilon|^2 dxds\le C
$$
with $\alpha\ge 0$,
then for any $r>0$,
\begin{equation}
\label{3-22}
\sup_{t\in(0, +\infty)}\|u_\varepsilon(\cdot,t)\|_{L^{r+1}}^{r+1} +
\sup_{t\in(0, +\infty)}\int_t^{t+1}\|\nabla u_\varepsilon^{\frac{m+r}2}(\cdot,s)\|_{L^{2}}^{2} ds
\le C_r,
\end{equation}
where $C_r$ depends on $r$, and is independent of $\varepsilon$.
\end{lemma}

{\it\bfseries Proof. } Noticing that $m-2\le 0$ and $\alpha\ge 0$,
we have
$$
|\nabla u_\varepsilon|^2\le (u_\varepsilon^{m-2}+u_\varepsilon^{\alpha})|\nabla u_\varepsilon|^2.
$$
By the assumption of this lemma, and recalling \eqref{3-2}, we obtain
\begin{equation}
\label{3-23}
\sup_{t\in(0,+\infty)}\int_t^{t+1}
\int_\Omega |\nabla u_\varepsilon|^2 dxds\le C.
\end{equation}
Applying $\nabla$ to the second equation in \eqref{3-1}, and
multiplying the resulting equation by $|\nabla v_\varepsilon|^{r-2}\nabla v_\varepsilon$
for any $r>2$, and using Lemma \ref{lem2-5}, we obtain
\begin{align}
&\frac1r\frac{d}{dt}\int_\Omega |\nabla v_\varepsilon|^{r} dx
+\int_{\Omega}|\nabla v_\varepsilon|^{r-2}|\nabla^2 v_\varepsilon|^2 dx+(r-2)
\int_{\Omega}|\nabla v_\varepsilon|^{r-2}(\nabla|\nabla v_\varepsilon|)^2 dx\nonumber
\\
&=-\int_{\Omega}|\nabla v_\varepsilon|^{r-2}\nabla v_\varepsilon\nabla(v_\varepsilon u_\varepsilon) dx+
\mu\int_{\Omega}|\nabla v_\varepsilon|^{r-2}\nabla v_\varepsilon\nabla(v_\varepsilon-v_\varepsilon^2) dx+
\int_{\partial\Omega}\frac{\partial(\nabla v_\varepsilon)}{\partial{\bf n}}|\nabla v_\varepsilon|^{r-2}
\nabla v_\varepsilon ds \nonumber
\\
&=-\int_{\Omega} u_\varepsilon|\nabla v_\varepsilon|^{r}dx-\int_{\Omega}v_\varepsilon
|\nabla v_\varepsilon|^{r-2}\nabla v_\varepsilon \nabla u_\varepsilon dx+\mu\int_{\Omega}(1-2v_\varepsilon)
|\nabla v_\varepsilon|^{r}dx
+\frac12
\int_{\partial\Omega}\frac{\partial(|\nabla v_\varepsilon|^2)}{\partial{\bf n}}
|\nabla v_\varepsilon|^{r-2}ds\nonumber
\\
\label{3-24}
&\le -\int_{\Omega} u_\varepsilon|\nabla v_\varepsilon|^{r}dx-\int_{\Omega}v_\varepsilon
|\nabla v_\varepsilon|^{r-2}\nabla v_\varepsilon \nabla u_\varepsilon dx+
\mu\int_{\Omega}|\nabla v_\varepsilon|^{r}dx+
\kappa\int_{\partial\Omega}|\nabla v_\varepsilon|^r ds.
\end{align}
By the boundary trace embedding inequalities, and combining with \eqref{3-2} and \eqref{3-5},
we conclude for any small $\eta>0$,
\begin{align*}
\kappa\int_{\partial\Omega}|\nabla v_\varepsilon|^r ds
\le \eta\|\nabla(|\nabla v_\varepsilon|^{\frac r2})\|_{L^2}^2
+C_\eta\||\nabla v_\varepsilon|^{\frac r2}\|_{L^{\frac4r}}^2
\\
\le \eta\|\nabla(|\nabla v_\varepsilon|^{\frac r2})\|_{L^2}^2
+\hat C_\eta.
\end{align*}
Substituting this inequality into \eqref{3-24} with $\eta$ appropriately small, we obtain
\begin{align}
&\frac1r\frac{d}{dt}\int_\Omega |\nabla v_\varepsilon|^{r} dx
+\int_{\Omega}|\nabla v_\varepsilon|^{r-2}|\nabla^2 v_\varepsilon|^2 dx+\frac{r-2}2
\int_{\Omega}|\nabla v_\varepsilon|^{r-2}(\nabla|\nabla v_\varepsilon|)^2 dx
+\int_{\Omega} u_\varepsilon|\nabla v_\varepsilon|^{r}dx\nonumber
\\
\label{3-25}
&\le \mu\int_{\Omega}|\nabla v_\varepsilon|^{r}dx-\int_{\Omega}v_\varepsilon
|\nabla v_\varepsilon|^{r-2}\nabla v_\varepsilon \nabla u_\varepsilon dx+C.
\end{align}
Taking $r=4$ in \eqref{3-25}, then for any small constant $\eta>0$, we have
\begin{align}
&\frac14\frac{d}{dt}\int_\Omega |\nabla v_\varepsilon|^{4} dx
+\int_{\Omega}|\nabla v_\varepsilon|^{2}|\nabla^2 v_\varepsilon|^2 dx+
\int_{\Omega}|\nabla v_\varepsilon|^{2}(\nabla|\nabla v_\varepsilon|)^2 dx
+\int_{\Omega} |\nabla v_\varepsilon|^{4}dx\nonumber
\\
&\le \eta\int_{\Omega}|\nabla v_\varepsilon|^6+C_\eta\int_{\Omega}|\nabla u_\varepsilon|^2 dx+
(\mu+1)\int_{\Omega} |\nabla v_\varepsilon|^{4}dx+C\nonumber
\\
\label{add3-26}
&\le 2\eta\int_{\Omega}|\nabla v_\varepsilon|^6+C_\eta\int_{\Omega}|\nabla u_\varepsilon|^2 dx+C.
\end{align}
Noticing that
\begin{align*}
\|\nabla v_\varepsilon\|_{r+2}^{r+2}
=\int_\Omega |\nabla v_\varepsilon|^{r}\nabla v_\varepsilon\nabla v_\varepsilon dx
=-\int_\Omega v_\varepsilon\left(|\nabla v_\varepsilon|^{r}\Delta v_\varepsilon+
r|\nabla v_\varepsilon|^{r-1}\nabla v_\varepsilon\nabla|\nabla v_\varepsilon| \right) dx
\\
\le C\left(\int_\Omega |\nabla v_\varepsilon|^{r-2}|\nabla^2 v_\varepsilon|^2 dx\right)^{\frac12}
\left(\int_\Omega |\nabla v_\varepsilon|^{r+2}dx\right)^{\frac12}+
\left(\int_{\Omega}|\nabla v_\varepsilon|^{r-2}(\nabla|\nabla v_\varepsilon|)^2 dx\right)^{\frac12}
\left(\int_\Omega |\nabla v_\varepsilon|^{r+2}dx\right)^{\frac12},
\end{align*}
then we have
\begin{equation}
\label{3-26}
\|\nabla v_\varepsilon\|_{r+2}^{r+2}\le
C\int_\Omega |\nabla v_\varepsilon|^{r-2}|\nabla^2 v_\varepsilon|^2 dx
+C\int_{\Omega}|\nabla v_\varepsilon|^{r-2}(\nabla|\nabla v_\varepsilon|)^2 dx.
\end{equation}
Substituting \eqref{3-26} into \eqref{add3-26} with $r=4$, and letting $\eta$ appropriately small, then
\begin{align*}
&\frac14\frac{d}{dt}\int_\Omega |\nabla v_\varepsilon|^{4} dx
+\frac12 \int_{\Omega}|\nabla v_\varepsilon|^{2}|\nabla^2 v_\varepsilon|^2 dx+
\frac12\int_{\Omega}|\nabla v_\varepsilon|^{2}(\nabla|\nabla v_\varepsilon|)^2 dx
+\int_{\Omega} |\nabla v_\varepsilon|^{4}dx
\\
&\le C_\eta\int_{\Omega}|\nabla u_\varepsilon|^2 dx+C.
\end{align*}
Using Lemma \ref{lem2-2} and \eqref{3-23}, we obtain
\begin{align}\label{3-27}
\sup_{t\in(0,+\infty)}\int_\Omega |\nabla v_\varepsilon|^{4} dx
\le C_1\sup_{t\in(0,+\infty)}\int_t^{t+1}\int_{\Omega}|\nabla u_\varepsilon|^2 dx+C_2
\le \tilde C.
\end{align}
By \eqref{3-19}, \eqref{3-27}, noticing that $r+2-m<m+r$, and taking $\eta$ appropriately small, we have
\begin{align*}
&\frac1{r+1}\frac{d}{dt}\int_\Omega u_\varepsilon^{r+1} dx
+\frac{rm}2\int_\Omega u_\varepsilon^{m+r-2}|\nabla u_\varepsilon|^2 dx+\int_\Omega u_\varepsilon^{r+1} dx
\\
\le & C\int_\Omega u_\varepsilon^{r+2-m}|\nabla v_\varepsilon|^{2}dx+C
\\
\le &C\|\nabla v_\varepsilon\|_{L^4}^2\|u_\varepsilon\|_{L^{2(r+2-m)}}^{r+2-m}+C
\\
\le &\eta\|u_\varepsilon\|_{L^{3(r+m)}}^{r+2-m}+C_\eta\|u_\varepsilon\|_{L^1}^{r+2-m}+C
\\
\le &\frac{rm}4\int_\Omega u_\varepsilon^{m+r-2}|\nabla u_\varepsilon|^2 dx+\hat C.
\end{align*}
Hence, \eqref{3-22} follows by using Lemma \ref{lem2-2}. \hfill $\Box$

\begin{lemma}
\label{lem3-7}
Assume (i) or (ii), and $m\le 2$.
Let $(u_\varepsilon, v_\varepsilon)$ be the solution of \eqref{3-1}. If
\begin{equation}\label{3-28}
\sup_{t\in(0,+\infty)}
\int_\Omega u_\varepsilon^{A_k(m-1)+1} dx+
\sup_{t\in(0,+\infty)}\int_t^{t+1}
\int_\Omega u_\varepsilon^{A_k(m-1)+m-2}|\nabla u_\varepsilon|^2 dxds\le C
\end{equation}
with $A_k(m-1)+m-2< 0$,
then there exists a constant $\tilde C$ independent of $\varepsilon$, such that
\begin{equation}\label{3-29}
\sup_{t\in(0,+\infty)}
\int_\Omega u_\varepsilon^{A_{k+1}(m-1)+1} dx+
\sup_{t\in(0,+\infty)}\int_t^{t+1}
\int_\Omega u_\varepsilon^{A_{k+1}(m-1)+m-2}|\nabla u_\varepsilon|^2 dxds\le \tilde C
\end{equation}
where $A_{k+1}=\frac23(m-1)A_k^2+(\frac{8m}{3}-2)A_k+2m-\frac13$.
\end{lemma}

{\it\bfseries Proof. } Recalling \eqref{3-25}, noticing that $2-m-A_k(m-1)>0$, and by Lemma \ref{lem2-1},
we obtain
\begin{align*}
&\frac1r\frac{d}{dt}\int_\Omega |\nabla v_\varepsilon|^{r} dx
+\int_{\Omega}|\nabla v_\varepsilon|^{r-2}|\nabla^2 v_\varepsilon|^2 dx+\frac{r-2}2
\int_{\Omega}|\nabla v_\varepsilon|^{r-2}(\nabla|\nabla v_\varepsilon|)^2 dx
+\int_{\Omega} u_\varepsilon|\nabla v_\varepsilon|^{r}dx
\\
&\le \mu\int_{\Omega}|\nabla v_\varepsilon|^{r}dx-\int_{\Omega}v_\varepsilon
|\nabla v_\varepsilon|^{r-2}\nabla v_\varepsilon \nabla u_\varepsilon dx+C
\\
&\le\mu\int_{\Omega}|\nabla v_\varepsilon|^{r}dx-\int_{\Omega}v_\varepsilon \left(u_\varepsilon^{\frac{A_k(m-1)+m-2}{2}}|\nabla u_\varepsilon|\right)
\left(u_\varepsilon^{\frac{2-m-A_k(m-1)}2}|\nabla v_\varepsilon|^{\frac{2-m-A_k(m-1)}{2}r}\right)
|\nabla v_\varepsilon|^{r-1-\frac{2-m-A_k(m-1)}{2}r} dx+C
\\
&\le \mu\int_{\Omega}|\nabla v_\varepsilon|^{r}dx+ C_\eta\int_\Omega u_\varepsilon^{A_k(m-1)+m-2}|\nabla u_\varepsilon|^2 dx
+\frac12\int_{\Omega} u_\varepsilon|\nabla v_\varepsilon|^{r}dx
\\
&+\eta\int_{\Omega}
|\nabla v_\varepsilon|^{\left(r-1-\frac{2-m-A_k(m-1)}{2}r\right)\frac{2}{(A_k+1)(m-1)}} dx+C
\end{align*}
for any sufficiently small constant $\eta>0$.
Noticing that $\left(r-1-\frac{2-m-A_k(m-1)}{2}r\right)\frac{2}{(A_k+1)(m-1)}=r+2$
when $r=2+2(A_k+1)(m-1)$, taking $\eta$ appropriately small in the above inequality,
and using \eqref{3-26}, we obtain
\begin{align*}
&\frac1{2+2(A_k+1)(m-1)}\frac{d}{dt}\int_\Omega |\nabla v_\varepsilon|^{2+2(A_k+1)(m-1)} dx
+\frac12\int_{\Omega}|\nabla v_\varepsilon|^{2(A_k+1)(m-1)}|\nabla^2 v_\varepsilon|^2 dx
\\
&+\frac{(A_k+1)(m-1)}2
\int_{\Omega}|\nabla v_\varepsilon|^{2(A_k+1)(m-1)}(\nabla|\nabla v_\varepsilon|)^2 dx
+\frac12\int_{\Omega} u_\varepsilon|\nabla v_\varepsilon|^{2+2(A_k+1)(m-1)}dx
\\
&\le  C\int_\Omega u_\varepsilon^{A_k(m-1)+m-2}|\nabla u_\varepsilon|^2 dx+C.
\end{align*}
Using \eqref{3-26}, \eqref{3-28} and Lemma \ref{lem2-2}, we obtain
\begin{align}
&\sup_{t\in(0,\infty)}\int_\Omega |\nabla v_\varepsilon|^{2+2(A_k+1)(m-1)} dx
+\sup_{t\in(0,\infty)}\int_t^{t+1}\int_{\Omega}|\nabla v_\varepsilon|^{2(A_k+1)(m-1)}|\nabla^2 v_\varepsilon|^2 dxds
\nonumber
\\
&+\sup_{t\in(0,\infty)}\int_t^{t+1}
\int_{\Omega}|\nabla v_\varepsilon|^{2(A_k+1)(m-1)}(\nabla|\nabla v_\varepsilon|)^2 dxds
+\sup_{t\in(0,\infty)}\int_t^{t+1}\int_{\Omega} u_\varepsilon|\nabla v_\varepsilon|^{2+2(A_k+1)(m-1)}dxds\nonumber
\\
\label{3-30}
&\le  C\sup_{t\in(0,\infty)}\int_t^{t+1}\int_\Omega u_\varepsilon^{A_k(m-1)+m-2}|\nabla u_\varepsilon|^2 dxds+C
\le C_k.
\end{align}
Recalling \eqref{3-19}, we see that for any small constant $\eta>0$,
\begin{align}
&\frac1{r+1}\frac{d}{dt}\int_\Omega u_\varepsilon^{r+1} dx
+\frac{2rm}{(r+m)^2}\int_\Omega|\nabla u_\varepsilon^{\frac{m+r}{2}}|^2 dx+\int_\Omega u_\varepsilon^{r+1} dx\nonumber
\\
\le & C\int_\Omega u_\varepsilon^{r+2-m}|\nabla v_\varepsilon|^{2}dx+C\nonumber
\\
\le & C\int_\Omega u_\varepsilon^{\frac1{1+(A_k+1)(m-1)}}|\nabla v_\varepsilon|^{2}
u_\varepsilon^{r+2-m-\frac1{1+(A_k+1)(m-1)}}dx+C\nonumber
\\
\label{3-31}
\le & C_\eta\int_{\Omega} u_\varepsilon|\nabla v_\varepsilon|^{2+2(A_k+1)(m-1)}dx+\eta\int_{\Omega}
u_\varepsilon^{\left(r+2-m-\frac1{1+(A_k+1)(m-1)}\right)\frac{1+(A_k+1)(m-1)}{(A_k+1)(m-1)}}+C.
\end{align}
By Gagliardo-Nirenberg interpolation inequality, for any $A>1$, it gives
\begin{align}
\|u_\varepsilon\|_{L^{r+m+\frac{2A}3}}^{r+m+\frac{2A}3}&
=\|u_\varepsilon^{\frac{m+r}2}\|_{L^{\frac{2(r+m+\frac{2A}3)}{m+r}}}^{\frac{2(r+m+\frac{2A}3)}{m+r}}
\le C_1\|u_\varepsilon^{\frac{m+r}{2}}\|_{L^{\frac{2A}{m+r}}}^{\frac{4A}{3(m+r)}}
\|\nabla u_\varepsilon^{\frac{m+r}{2}}\|_{L^2}^2+C_2\|u_\varepsilon\|_{L^1}^{r+m+\frac{2A}3}\nonumber
\\
\label{3-32}
&\le C_1 \|u_\varepsilon\|_{L^A}^{\frac{2A}{3}}\|\nabla u_\varepsilon^{\frac{m+r}{2}}\|_{L^2}^2+C_3.
\end{align}
By \eqref{3-28} and \eqref{3-32}, we get
\begin{align}\label{3-33}
\int_\Omega u_\varepsilon^{r+m+\frac{2A_k(m-1)+2}3} dx\le
C_4\int_\Omega|\nabla u_\varepsilon^{\frac{m+r}{2}}|^2 dx+C_5.
\end{align}
Combining \eqref{3-31} and \eqref{3-33}, and noticing that
$$
\left(r+2-m-\frac1{1+(A_k+1)(m-1)}\right)\frac{1+(A_k+1)(m-1)}{(A_k+1)(m-1)}=r+m+\frac{2A_k(m-1)+2}3
$$
if $r=A_{k+1}(m-1)$ with $A_{k+1}=\frac23(m-1)A_k^2+(\frac{8m}{3}-2)A_k+2m-\frac13$,
and \eqref{3-29} is arrived by lemma \ref{lem2-2}. \hfill $\Box$

\medskip

Next, we pay our attention to  the  convergence and divergence properties of the sequence $\{A_k\}$.

\begin{lemma}
\label{lem3-8}
Assume that $1<m<2$.
Let $A_{k+1}=\frac23(m-1)A_k^2+(\frac{8m}{3}-2)A_k+2m-\frac13$ with $A_1=1$, then we have

i) when $m>\frac{11}4-\sqrt 3$, $\lim_{k\to\infty}A_k=+\infty$;

ii) when $1<m\le \frac{11}4-\sqrt 3$, $\lim_{k\to\infty}A_k=A^*$,
where
$A^*=\frac{9-8m-\sqrt{16m^2-88m+73}}{4(m-1)}>0$.
\end{lemma}

{\it\bfseries Proof.}
Firstly, by a direct calculation, we see that $A_2=\frac{16m}3-3>A_1$.
If $A_{k}>A_{k-1}$,
it is easy to see that $A_{k+1}=\frac23(m-1)A_k^2+(\frac{8m}{3}-2)A_k+2m-\frac13
>\frac23(m-1)A_{k-1}^2+(\frac{8m}{3}-2)A_{k-1}+2m-\frac13=A_k$, which implies that $\{A_k\}$
is an increasing sequence.

In what follows, we first show i). Suppose the contrary, then by  monotone convergence theorem, there
exists a positive constant $A$ such that
$$
\lim_{k\to\infty}A_k=A.
$$
Then we have
$$
A=\frac23(m-1)A^2+(\frac{8m}{3}-2)A+2m-\frac13,
$$
that is
\begin{equation}\label{3-34}
\frac23(m-1)A^2+(\frac{8m}{3}-3)A+2m-\frac13=0.
\end{equation}
By  Weda's Theorem, let
$$
\Delta=(\frac{8m}{3}-3)^2-\frac83(m-1)(2m-\frac13)=\frac{16m^2-88m+73}{9}
$$
it is not difficult to verify that when $\frac{11}4-\sqrt 3<m<\frac{11}4+\sqrt 3$, $\Delta<0$,
that is \eqref{3-34} has no solution, which is a contradiction. Then i) is proved.

Next, we show ii). We first verify that $A_k$ has an upper bound.
By a direct calculation, we see that $m\le \frac{11}4-\sqrt 3<\frac{13}{12}$,
then it is not difficult to verify that $A*>1=A_1$. By recurrence method, we see that
if $A_k<A^*$, then
$$
A_{k+1}=\frac23(m-1)A_k^2+(\frac{8m}{3}-2)A_k+2m-\frac13<\frac23(m-1)A^{*2}+(\frac{8m}{3}-2)A^*+2m-\frac13=A^*,
$$
which means that $A^*$ is an upper bound of $\{A_k\}$, then by  monotone convergence theorem,
there exists a positive constant $A\le A^*$ such that
$$
\lim_{k\to\infty}A_k=A,
$$
and $A$ satisfies the equation \eqref{3-34}, and ii) is proved. \hfill $\Box$

\medskip

Combining Lemma \ref{lem3-5}--Lemma \ref{lem3-8}, we conclude the following lemma.

\begin{lemma}
\label{lem3-9}
Assume  (i) or (ii), and
let $(u_\varepsilon, v_\varepsilon)$ be the solution of \eqref{3-1}. If $\frac{11}4-\sqrt 3<m\le 2$,
then for any $r>0$, we have
\begin{equation}
\label{3-35}
\sup_{t\in(0, +\infty)}\|u_\varepsilon(\cdot,t)\|_{L^{r+1}}^{r+1} +
\sup_{t\in(0, +\infty)}\int_t^{t+1}\|\nabla u_\varepsilon^{\frac{m+r}2}(\cdot,s)\|_{L^{2}}^{2} ds
\le C_r
\end{equation}
and
\begin{equation}
\label{3-36}
\sup_{t\in(0, +\infty)}\|v_\varepsilon\|_{W_r^{2,1}Q_1(t)}\le \tilde C_r,
\end{equation}
where $C_r$, $\tilde C_r$ depend on $r$, and are independent of $\varepsilon$.
\end{lemma}

{\it\bfseries Proof. } By Lemma \ref{lem3-5} and \eqref{lem3-7}, we consider the  sequence of iterations $\{A_k\}$
with $A_1=1$. By Lemma \ref{lem3-8}, we have $A_k\to +\infty$ as $k\to +\infty$, if $m>\frac{11}4-\sqrt 3$.
Then there exists $K>0$ such that $A_K(m-1)+m-2\ge 0$, thus by Lemma \eqref{lem3-6}, \eqref{3-35} is arrived.
And \eqref{3-36} is subsequently obtained by using Lemma \ref{lem2-3}.
\hfill $\Box$

\begin{remark}
By Lemma \ref{lem3-5}--Lemma \ref{lem3-8}, we also see that
 if $m\le \frac{11}4-\sqrt 3$. Then for any $A<A^*$, we have
\begin{equation}\label{3-37}
\sup_{t\in(0,+\infty)}
\int_\Omega u_\varepsilon^{A(m-1)+1} dx+
\sup_{t\in(0,+\infty)}\int_t^{t+1}
\int_\Omega u_\varepsilon^{A(m-1)+m-2}|\nabla u_\varepsilon|^2 dxds\le C,
\end{equation}
where $A^*$ is defined in ii) of Lemma \ref{lem3-8}. It is not difficulty to verify that $A^*>5$ for any
$1<m\le \frac{11}4-\sqrt 3$. Thus \eqref{3-37} holds for $A=5$. But it is hard to go further
for regularity estimate
of $u_\varepsilon$.
\end{remark}

Using these lemmas, we further have
\begin{lemma}
\label{lem3-10}
Assume  (i) or (ii), and
let $(u_\varepsilon, v_\varepsilon)$ be the solution of \eqref{3-1}. Then
when $m>\frac{11}4-\sqrt 3$,
\begin{align}
\label{3-38}
&\sup_{t\in(0,+\infty)}\|v_\varepsilon(\cdot, t)\|_{W^{1,\infty}}\le C,
\\
\label{3-39}
& \sup_{t\in(0,+\infty)}\|u_\varepsilon(\cdot, t)\|_{L^{\infty}}\le C,
\end{align}
where the constants $C$ are independent of $\varepsilon$.
\end{lemma}

{\it\bfseries Proof. }  By $t$-anisotropic embedding theorem,
that is $W_p^{2,1}\hookrightarrow C^{\alpha, \frac{\alpha}2}$
for any $0<\alpha\le 2-\frac{5}{p}$, we obtain \eqref{3-38} by Lemma \ref{lem3-4} and Lemma\ref{lem3-9}. By \eqref{3-38},
and similarly to the proof of in \cite{J2},  we can obtain the $L^\infty$ estimate
of $u_\varepsilon$ by using Moser iteration method.  \hfill $\Box$

\begin{lemma}
\label{lem3-11}
Assume  (i) or (ii),   and let
$(u_\varepsilon, v_\varepsilon)$ be the solution of \eqref{3-1}. Then when $m>\frac{11}4-\sqrt 3$,
\begin{align}
\label{3-40}
\sup_{t\in(0, +\infty)}\int_\Omega|\nabla u_\varepsilon^m|^2 dx+
\varepsilon\sup_{t\in(0, +\infty)}\int_t^{t+1}\int_\Omega \left|\frac{\partial u_\varepsilon}{\partial t}\right|^2 dxds
+\sup_{t\in(0, +\infty)}\int_t^{t+1}\int_\Omega u_\varepsilon^{m-1}
\left|\frac{\partial u_\varepsilon}{\partial t}\right|^2 dxds\le C,
\end{align}
where $C$ is independent of $\varepsilon$.
\end{lemma}

{\it\bfseries Proof.} Multiplying the first equation of \eqref{3-1} by
$\frac{\partial(\varepsilon u_\varepsilon+u_\varepsilon^m)}{\partial t}$,
and integrating it over $\Omega$ gives
\allowdisplaybreaks
\begin{align*}
&\frac1{2}\frac{d}{dt}\int_\Omega|\nabla (\varepsilon u_\varepsilon+u_\varepsilon^m)|^2 dx+
\varepsilon \int_\Omega \left|\frac{\partial u_\varepsilon}{\partial t}\right|^2 dx
+m\int_\Omega u_\varepsilon^{m-1} \left|\frac{\partial u_\varepsilon}{\partial t}\right|^2 dx
\\
&\le - \chi\int_\Omega \nabla\cdot(u_{\varepsilon}\nabla v_{\varepsilon})
\frac{\partial(\varepsilon u_\varepsilon+u_\varepsilon^m)}{\partial t}dx
+\int_\Omega (\xi u_{\varepsilon}v_{\varepsilon}-\rho u_{\varepsilon}-\varepsilon u_{\varepsilon}^2)
\frac{\partial(\varepsilon u_\varepsilon+u_\varepsilon^m)}{\partial t}dx
\\
&\le  m\chi^2\int_\Omega|\nabla\cdot(u_{\varepsilon}\nabla v_{\varepsilon})|^2
u_\varepsilon ^{m-1} dx +\varepsilon\chi^2\int_\Omega|\nabla\cdot(u_{\varepsilon}\nabla v_{\varepsilon})|^2dx
+\frac{\varepsilon}2 \int_\Omega \left|\frac{\partial u_\varepsilon}{\partial t}\right|^2 dx
\\
&+\varepsilon\int_\Omega (\xi u_{\varepsilon}v_{\varepsilon}-\rho u_{\varepsilon}-\varepsilon u_{\varepsilon}^2)^2dx
+m\int_\Omega u_{\varepsilon}^{m-1}(\xi u_{\varepsilon}v_{\varepsilon}-\rho u_{\varepsilon}-\varepsilon u_{\varepsilon}^2)^2 dx
+\frac{m}2\int_\Omega u_\varepsilon^{m-1} \left|\frac{\partial u_\varepsilon}{\partial t}\right|^2 dx,
\end{align*}
recalling \eqref{3-38} and \eqref{3-39}, we further have
\begin{align*}
&\frac1{2}\frac{d}{dt}\int_\Omega|\nabla (\varepsilon u_\varepsilon+u_\varepsilon^m)|^2 dx+
\frac{\varepsilon}2 \int_\Omega \left|\frac{\partial u_\varepsilon}{\partial t}\right|^2 dx
+\frac{m}2\int_\Omega u_\varepsilon^{m-1} \left|\frac{\partial u_\varepsilon}{\partial t}\right|^2 dx
+\int_\Omega|\nabla (\varepsilon u_\varepsilon+u_\varepsilon^m)|^2 dx
\\
&\le  m\chi^2\int_\Omega|\nabla\cdot(u_{\varepsilon}\nabla v_{\varepsilon})|^2
u_\varepsilon ^{m-1} dx +\varepsilon\chi^2\int_\Omega|\nabla\cdot(u_{\varepsilon}\nabla v_{\varepsilon})|^2dx
+\int_\Omega|\nabla (\varepsilon u_\varepsilon+u_\varepsilon^m)|^2 dx+C
\\
&\le  C\int_\Omega (|\nabla u_\varepsilon^{\frac m2}|^2+|\Delta v_{\varepsilon}|^2)dx+
C\varepsilon\int_\Omega \frac{|\nabla u_\varepsilon|^2}{u_\varepsilon} dx+C.
\end{align*}
Recalling \eqref{3-5},  \eqref{3-17}, \eqref{3-36}, \eqref{3-39}  and using Lemma \ref{lem2-2}, we obtain
$$
\sup_{t\in(0, +\infty)}\int_\Omega|\nabla (\varepsilon u_\varepsilon+u_\varepsilon^m)|^2 dx+
\varepsilon\sup_{t\in(0, +\infty)}\int_t^{t+1}\int_\Omega \left|\frac{\partial u_\varepsilon}{\partial t}\right|^2 dxds
+\sup_{t\in(0, +\infty)}\int_t^{t+1}\int_\Omega u_\varepsilon^{m-1}
\left|\frac{\partial u_\varepsilon}{\partial t}\right|^2 dxds\le C.
$$
and \eqref{3-40} is obtained. \hfill $\Box$

\medskip

Completely similar to the proof as Lemma \ref{lem3-3}-Lemma \ref{lem3-11}, for the case
$\xi\mu>0$ and $\rho=0$, we also have
\begin{lemma}
\label{lem3-12}
Assume (iii) $\xi\mu>0$ and $\rho=0$,  and let
$(u_\varepsilon, v_\varepsilon)$ be the solution of \eqref{3-1}. Then
for any $T>0$,
\begin{align}
\label{3-41}
\sup_{t\in(0,T)}\left(\|u_\varepsilon(\cdot, t)\|_{L^{\infty}}+\|\nabla u_\varepsilon^m\|_{L^2}+
\|v_\varepsilon(\cdot, t)\|_{W^{1,\infty}}\right)
+\iint_{Q_T}\left(\varepsilon|\nabla u_{\varepsilon}|^2+|\nabla u_{\varepsilon}^{\frac m2}|^2
+ \left|\frac{\partial u_\varepsilon^{\frac {m+1}2}}{\partial t}\right|^2\right) dxdt\le C,
\end{align}
and
\begin{equation}
\label{3-42}
\|v_\varepsilon\|_{W_p^{2,1}(Q_T)}\le C_p, \  \text{for any}\   p\in (1, +\infty),
\end{equation}
where  the constants $C_p$ and $C$ depend on  $T$,  $C_p$ depends on $p$, and both of them are independent
of $\varepsilon$.
\end{lemma}

{\it \bfseries Proof of Theorem \ref{thm-1} and Theorem \ref{thm-2}.}
Since $(u_\varepsilon, v _\varepsilon, \omega_\varepsilon)$ is the classical solution of
\eqref{3-1}, then we have
\begin{align*}
&-\iint_{Q_T}u_\varepsilon\varphi_t dxdt-\int_\Omega u_\varepsilon(x,0)\varphi(x,0)dx
+ \varepsilon \iint_{Q_T}\nabla u_\varepsilon\nabla\varphi
+\iint_{Q_T} \nabla u_\varepsilon^m\nabla\varphi dxdt
\\
&=\chi\iint_{Q_T}u_\varepsilon\nabla v_\varepsilon\nabla\varphi dxdt+
\xi \iint_{Q_T} u_{\varepsilon}v_{\varepsilon}\varphi dxdt
-\iint_{Q_T}(\rho u_{\varepsilon}+\varepsilon u_{\varepsilon}^2)\varphi dxdt
\\
&-\iint_{Q_T}v_\varepsilon\varphi_tdxdt-\int_\Omega v_\varepsilon(x,0)\varphi(x,0)dx
+\iint_{Q_T}\nabla v_\varepsilon\nabla\varphi dxdt=\iint_{Q_T}
(-v_{\varepsilon}u_{\varepsilon}
+\mu v_{\varepsilon}(1-v_{\varepsilon}))\varphi dxdt,
\end{align*}
for any $\varphi\in C^\infty(\overline Q_T)$ with  $\varphi(x,T)=0$.
By \eqref{3-5}, \eqref{3-39}, or \eqref{3-41}, we also have
$$
\sqrt \varepsilon\nabla u_\varepsilon\in L^2(Q_T),
$$
it implies that  when $\varepsilon\to 0$,
$$
\varepsilon \nabla u_\varepsilon \to 0, \ \text{in  $L^2(Q_T)$}.
$$
By Sobolev compact embedding theorem, and using Lemma \ref{lem3-3},
Lemma \ref{lem3-4}, Lemma \ref{lem3-9},
Lemma \ref{lem3-10}, Lemma \ref{lem3-11} and Lemma \ref{lem3-12},
for any $T>0$, letting $\varepsilon\to 0$, we have
\begin{align*}
&u_\varepsilon \to u, \varepsilon u_\varepsilon^2 \to 0, \ \text{in $L^p(Q_T)$, for any $p\in (1,+\infty)$},
\\
&u_\varepsilon \stackrel{*}{\rightharpoonup} u, \ \text{in $L^\infty(Q_T)$},
\\
&\nabla u_\varepsilon^m \to \nabla u^m, \ \text{in  $L^2(Q_T)$},
\\
&v_\varepsilon \to v, \ \text{uniformly},
\\
& v_\varepsilon \rightharpoonup v, \ \text{in $W^{2,1}_p(Q_T)$, for any $p\in (1,+\infty)$},
\end{align*}
which means $(u, v)$ is the global weak solution of \eqref{1-1} with \eqref{1-4}--\eqref{1-9} hold. \hfill $\Box$

\section{Large Time Behavior of Solutions}
To investigate the large time behavior of solutions of the problem \eqref{1-1}, we need the following two lemmas \cite{J3}.
\begin{lemma}
\label{lem4-1}
Assume that $f\ge 0$, $f(t)\in L^1(T,\infty)$ for some constant $T>0$, and
$$
f(t)-f(s)\le A(t-s) \quad (\ \text{or}\ \ge -A(t-s)),\quad \text{for any }\ t>s>T.
$$
Then
$$
\lim_{t\to\infty}f(t)=0.
$$
\end{lemma}
\begin{lemma}
\label{lem4-2}
Assume that $f(t), g(t)\ge0$, $\lim_{t\to\infty}g(t)=0$, and $f(t)\in L^1(T,+\infty)$ with some constant $T\ge 0$.
Let $F(t)=f(t)-g(t)$, and that
$$
F(t)-F(s)\ge -A(t-s), \forall t>s>T,
$$
then
$$
\lim_{t\to\infty}f(t)=0.
$$
\end{lemma}

Firstly for the case $\mu=0$, we have the following lemma.
\begin{lemma}
\label{lem4-3}
Assume that $\mu=0$. Let $(u, v)$ with $u\in\mathcal X_1, v\in \mathcal X_2$  be the global solution.
Then we have
\begin{equation}
\label{4-1}
\int_0^\infty\int_\Omega(|\Delta v|^2+|\nabla v|^2+uv+u^{m+r-2}|\nabla u|^2+\rho u)dxdt\le C_r, \quad
\text{for any $r>(m-2)_+$},
\end{equation}
where $C_r$ is a constant depends on $r$ and the initial datum $(u_0, v_0)$.
In particular,
\begin{equation}
\label{4-2}
\lim_{t\to\infty}\int_\Omega(|\nabla v|^2+\rho u) dx=0.
\end{equation}
\end{lemma}

{\it\bfseries Proof.}
Integrating the second equation of \eqref{1-1} over $\Omega$, we see that
$$
\frac{d}{dt}\int_{\Omega} vdx+\int_{\Omega}uvdx=0,
$$
integrating this equality from $0$ to $\infty$, we obtain
\begin{equation}\label{4-3}
\int_0^\infty\int_\Omega uv dxdt\le \int_{\Omega} v_0dx.
\end{equation}
Multiplying the second equation of \eqref{1-1} by $v$, we have
\begin{align*}
\frac{1}{2}\frac{d}{dt}\int_{\Omega}v^2dx+\int_\Omega|\nabla v|^2dx+\int_{\Omega}v^2u dx=0,
\end{align*}
integrating this equality from $0$ to $\infty$, we obtain
\begin{align}
\label{4-4}
\int_0^\infty\int_\Omega(|\nabla v|^2+v^2u)dxdt\le \int_\Omega |v_0|^2.
\end{align}
Multiplying the second equation of \eqref{1-1} by  $\Delta v$, and using Young's inequality, we obtain
\begin{align}
\label{4-5}
\frac{1}{2}\frac{d}{dt}\int_{\Omega}|\nabla v|^2dx+\frac12\int_\Omega|\Delta v|^2dx\le \frac12\int_\Omega v^2u^2dx
\le \frac12\|u\|_{L^\infty}\int_\Omega v^2udx,
\end{align}
integrating this equality from $0$ to $\infty$, we obtain
\begin{align}
\label{4-6}
\int_0^\infty\int_\Omega|\Delta v|^2dxdt\le \int_\Omega |\nabla v_0|^2+\sup_{t}\|u\|_{L^\infty}\int_0^\infty\int_\Omega v^2udxdt.
\end{align}
Integrating the first equation  of \eqref{1-1} over $\Omega$, we obtain
\begin{equation}
\label{4-7}
\frac{d}{dt}\int_{\Omega} udx+\rho\int_\Omega udx=\xi\int_{\Omega}uvdx,
\end{equation}
which implies
\begin{equation}
\label{4-8}
\int_{\Omega} u(x,t)dx+\rho\int_0^\infty\int_{\Omega}udxdt=\xi\int_0^\infty\int_{\Omega}uvdxdt+\int_{\Omega} u_0dx\le C.
\end{equation}
Multiplying the first equation of \eqref{1-1} by $u^r$ for any $r>(m-2)_+$, we obtain
\begin{align*}
&\frac{1}{r+1}\frac{d}{dt}\int_{\Omega}u^{r+1}dx+mr\int_\Omega u^{m+r-2}|\nabla u|^2dx+
\rho\int_\Omega u^{r+1}dx
=\chi r\int_\Omega u^r\nabla u\nabla vdx+\xi\int_\Omega u^{r+1}vdx
\\
&\le \frac{mr}2\int_\Omega u^{m+r-2}|\nabla u|^2dx+C\int_\Omega u^{r+2-m}|\nabla v|^2dx+\xi\int_\Omega u^{r+1}vdx,
\\
&\le\frac{mr}2\int_\Omega u^{m+r-2}|\nabla u|^2dx+C\int_\Omega |\nabla v|^2dx
+C\int_\Omega uvdx,
\end{align*}
by \eqref{4-3} and \eqref{4-4} and the boundedness of $u$ and $v$, we obtain
\begin{equation}
\label{4-9}
\int_0^\infty\int_\Omega u^{m+r-2}|\nabla u|^2dxdt\le C.
\end{equation}
Combining \eqref{4-3},\eqref{4-4}, \eqref{4-6}. \eqref{4-8} and \eqref{4-9},  \eqref{4-1} is proved.

On the other hands, by \eqref{4-5} and \eqref{4-7}, and using the boundedness of $u$ and $v$, we also have
$$
\frac{d}{dt}\int_{\Omega}(|\nabla v|^2+\rho u)dx\le C,
$$
which implies that
$$
\int_{\Omega}(|\nabla v(x,t)|^2+\rho u(x,t))dx
-\int_{\Omega}(|\nabla v(x,s)|^2+\rho u(x,s))dx\le C(t-s).
$$
Hence, by Lemma \ref{lem4-1}  and the inequality \eqref{4-1}, we obtain \eqref{4-2}.
The proof is complete. \hfill $\Box$

\begin{lemma}
\label{lem4-4}
Assume that $\mu=0$, $\rho=0$, $\xi\ge 0$, $u_0\not\equiv0$. Let $(u, v)$ with $u\in\mathcal X_1, v\in
\mathcal X_2$  be
the global solution.
Then we have
\begin{align}
\label{4-10}
\lim_{t\to\infty}\|v\|_{L^\infty}=0, \qquad \lim_{t\to\infty}\|u-A\|_{L^p}=0, \quad\text{for any $p>1$}.
\end{align}
where $A=\frac{1}{|\Omega|}\int_\Omega(u_0+\xi v_0)dx>0$.
\end{lemma}

{\it\bfseries Proof.}
We denote
$$
a(t)=\frac{1}{|\Omega|}\int_\Omega u(x,t)dx, \qquad b(t)=\frac{1}{|\Omega|}\int_\Omega v(x,t)dx.
$$
It is easy to see that
$$
a'(t)=\frac{\xi}{|\Omega|}\int_\Omega uvdx\ge 0,
\qquad b'(t)=-\frac{1}{|\Omega|}\int_\Omega uvdx\le 0,
$$
$$
(a+\xi b)'(t)=0,
$$
which means that
$$
a(t)+\xi b(t)\equiv A.
$$
Note that $a(t)$ is monotonically increasing and bounded above, $b(t)$ is monotonically
decreasing and bounded below, then there exists two constants $a^*>0$ and $b^*\ge 0$ such that
$a^*+b^*=A$, and
$$
\lim_{t\to\infty}a(t)=a^*, \qquad \lim_{t\to\infty} b(t)=b^*.
$$
By Poincar\'e inequality, we have
$$
\|v-b(t)\|_{L^2}\le C\|\nabla v\|_{L^2},
$$
and note that
$$
\|v-b^*\|_{L^2}\le \|v-b(t)\|_{L^2}+|\Omega|^{1/2}|b(t)-b^*|,
$$
combining with \eqref{4-2}, we have
$$
\lim_{t\to\infty}\|v-b^*\|_{L^2}=0.
$$
and we further have
\begin{equation}\label{4-11}
\lim_{t\to\infty}\|v-b^*\|_{L^\infty}=0
\end{equation}
since
$$
\|v-b^*\|_{L^\infty}\le C_1\|v-b^*\|_{L^2}^{\frac25}\|\nabla v\|_{L^\infty}^{\frac35}
+C_2\|v-b^*\|_{L^2}.
$$
Next, we show $b^*=0$. Suppose to the contrary, that is $b^*>0$.
By \eqref{1-6}, we also have $v\in C^{\alpha, \frac{\alpha}2}(Q)$ for some positive constant
$\alpha>0$, then there exists $T_0>0$, such that
$$
v(x,t)>\frac{b^*}2, \qquad \text{for any $t>T_0$}.
$$
Then we have
\begin{align*}
b'(t)==-\frac{1}{|\Omega|}\int_\Omega uvdx\le -\frac{b^*}{2|\Omega|}\int_\Omega udx=-\frac{b^*}2a(t)\le
-\frac{a^*b^*}2<0,
\end{align*}
which implies that $b(t)\to 0$, that is $b^*=0$. It is a contradiction.
The first equality of \eqref{4-10} is proved. It also implies that $a^*=A$,
that is
\begin{equation}
\label{4-12}
\lim_{t\to\infty}a(t)=A.
\end{equation}
Next, we show the second equality of \eqref{4-10}. Denote
$$
\alpha^m(t)=\frac{1}{|\Omega|}\int_\Omega u^m(x,t)dx,
$$
then we have
\begin{equation}
\label{4-13}
\alpha^m(t)=\frac1{|\Omega|}\int_\Omega u^{m}dx\ge \left(\frac1{|\Omega|}\int_\Omega udx\right)^{m}=a(t)^m.
\end{equation}
By \eqref{4-1} with $r=m$ and Poincar\'e inequality, we also have
\begin{equation}
\label{4-14}
\int_0^\infty\int_\Omega |u^{m}-\alpha^m|^2dxdt\le \int_0^\infty\int_\Omega |\nabla u^{m}|^2dxdt\le C.
\end{equation}
By a direct calculation, and using \eqref{1-4}, \eqref{1-5}, we obtain
\begin{align*}
&\frac{d}{dt}\int_\Omega|u^m-\alpha^m|^2 dx=\frac{d}{dt}\int_\Omega u^{2m}dx-2|\Omega|\alpha^m(\alpha^m)'
\\
&=2m^2(m-1)\alpha^m\int_\Omega u^{2m-3}|\nabla u|^2dx-2(2m-1)\int_\Omega u^{m-1}|\nabla u^m|^2dx
-(2m-1)\chi\int_\Omega u^{2m}\Delta vdx
\\
&+2(m-1)\chi \alpha^m\int_\Omega u^{m}\Delta vdx
+2 m\xi\int_\Omega u^{2m}vdx-2m\alpha^m\xi\int_\Omega u^{m}vdx-2m\rho\int_\Omega u^{2m}dx+2m\rho\alpha^m\int_\Omega u^{m}dx
\\
&\ge -C-\int_\Omega|\Delta v|^2dx.
\end{align*}
Combining with \eqref{4-5}, we obtain
\begin{align*}
\frac{d}{dt}\int_\Omega(|u^m-\alpha^m|^2-|\nabla v|^2) dx\ge -C,
\end{align*}
which implies that for any $t>s>0$,
$$
\int_\Omega(|u^m(x,t)-\alpha^m(t)|^2-|\nabla v(x,t)|^2) dx-\int_\Omega(|u^m(x,s)-\alpha^m(s)|^2-|\nabla v(x,s)|^2) dx\ge -C(t-s).
$$
Using \eqref{4-2}, \eqref{4-14} and Lemma \ref{lem4-2}, we obtain
\begin{equation}
\label{4-15}
\lim_{t\to\infty}\int_\Omega |u^m(x,t)-\alpha^m(t)|^2dx=0.
\end{equation}
Noticing that $\frac{u^m-\alpha^m}{u-\alpha}\ge \alpha^{m-1}$,  and using \eqref{4-13},  we have
\begin{equation}
\label{4-16}
a(t)^{2m-2}\int_\Omega |u-\alpha(t)|^2 dx\le \alpha^{2m-2}(t)\int_\Omega |u-\alpha(t)|^2 dx \le \int_\Omega |u^m-\alpha(t)^m|^2 dx.
\end{equation}
Combining \eqref{4-12}, \eqref{4-15} and \eqref{4-16}, we obtain
\begin{equation}
\label{4-17}
|a(t)-\alpha(t)|^2=\left|\frac 1{|\Omega|}\int_\Omega (u-\alpha) dx\right|^2\le
\int_\Omega |u-\alpha|^2 dx\to 0, \quad \text{as} \quad t\to\infty.
\end{equation}
Therefore, we have
\begin{equation}
\label{4-18}
\|u-A\|_{L^2}\le \|u-\alpha(t)\|_{L^2}+|\Omega|^{\frac12}|\alpha(t)-a(t)|+|\Omega|^{\frac12}|a(t)-A|\to 0,\quad \text{as} \quad t\to\infty.
\end{equation}
If $p<2$,  clearly, we have
$$
\|u-A\|_{L^p}\le \|u-A\|_{L^2};
$$
for any $p>2$, we see that
$$
\|u-A\|_{L^p}\le \|u-A\|_{L^2}^{\frac{2}p}\|u-A\|_{L^\infty}^{\frac{p-2}p}\le C\|u-A\|_{L^2}^{\frac{2}p},
$$
and the proof is complete. \hfill $\Box$
\begin{lemma}
\label{lem4-5}
Assume that $\mu=0$, $\rho>0$, $\xi\ge 0$, and  $v_0\not\equiv 0$. Let $(u, v)$ with $u\in\mathcal X_1, v\in
\mathcal X_2$  be
the global solution.
Then there exists a constant $B\in \left(0,\ \frac{1}{|\Omega|}\int_\Omega v_0 dx\right)$, such that
\begin{align}
\label{4-19}
\lim_{t\to\infty}\|v-B\|_{L^\infty}=0,  \qquad \lim_{t\to\infty}\|u\|_{L^p}=0 \quad \text{for any $p>1$}.
\end{align}
\end{lemma}

{\it\bfseries Proof.}
Note that $b(t)$ is monotonically
decreasing and bounded below, where $b(t)$ is defined as the proof of Lemma \ref{lem4-3}. Then there exists $B\ge 0$
such that
$$
\lim_{t\to\infty}b(t)=B.
$$
Similar to the proof of Lemma \ref{lem4-4}, we have
\begin{equation}
\label{add4-20}
\lim_{t\to\infty}\|v-B\|_{L^\infty}=0.
\end{equation}
In what follows, we show that $B>0$.
Let's consider the following problem
\begin{align}
\label{4-20}\left\{
\begin{aligned}
&\tilde v_t-\Delta \tilde v+M\tilde v=0,
\\
&\frac{\partial{\tilde v}}{\partial \bf n}\Big|_{\partial\Omega}=0,
\\
&\tilde v(x,0)=v_0(x).
\end{aligned}\right.
\end{align}
By \cite{HPW}, there exists $\Gamma_0>0$, such  that
$$
\tilde v(x,t)=e^{-Mt}e^{t\Delta}v_0\ge e^{-Mt}\Gamma_0\int_{\Omega}v_0(x) dx, \ \text{for  any}\ t\ge 1,
$$
where $\{e^{t\Delta}\}_{t\ge 0}$ is the Neumann heat semigroup.
Taking $M\ge \sup_{t}\|u(\cdot,t)\|_{L^\infty}$,  and by comparison, we have
\begin{equation}
\label{4-21}
v(x,t)\ge \tilde v(x,t)\ge e^{-Mt}\Gamma_0\int_{\Omega}v_0(x) dx, \ \text{for  any}\ t\ge 1.
\end{equation}
Multiplying the second equation of \eqref{1-1} by $\frac1{v}$,  integrating it over $\Omega\times (1,t)$,
and combining with \eqref{4-1}and \eqref{4-21}, we obtain
\begin{align}
\label{4-22}
\int_{\Omega}\ln v(x,t) dx=\int_{\Omega}\ln v(x,1) dx+\int_1^{t}\int_{\Omega}\frac{|\nabla v|^2}{v^2}dxds-
\int_1^{t}\int_{\Omega} u dxds\ge -C
\end{align}
for some positive constant $C$. Recalling \eqref{add4-20},  it implies that $B>0$.

Next, we show the second  limit equality.
By \eqref{4-2}, we see that
$$
\lim_{t\to\infty}\|u\|_{L^1}=0
$$
since $\rho>0$.
Thus, we further have
$$
\|u\|_{L^p}\le \|u\|_{L^1}^{\frac1p} \|u\|_{L^\infty}^{\frac{p-1}p}\to 0, \ \text{as $t\to\infty$}
$$
for any $p>1$, and this lemma is proved. \hfill $\Box$

\medskip

Next, we turn our attention to discuss the asymptotic behavior of solutions to the problem \eqref{1-1} in the case of $\mu>0$.

\begin{lemma}
\label{lem4-6}
Assume that $\mu>0$, $\rho>0$, $0\le \xi< \rho$, and  $v_0\not\equiv 0$. Let $(u, v)$ with $u\in\mathcal X_1, v\in
\mathcal X_2$  be
the global solution.
Then
\begin{align}
\label{4-23}
\lim_{t\to\infty}\|v-1\|_{L^\infty}=0,  \qquad \lim_{t\to\infty}\|u\|_{L^p}=0 \quad \text{for any $p>1$}.
\end{align}
\end{lemma}

{\it\bfseries Proof.}  Firstly, similar to the proof of \eqref{4-21}, there exists $\delta>0$ such that
$$
v(x,1)>\delta,
$$
for any $x\in\Omega$ since $\int_{\Omega}v_0 dx>0$.  Let
$$
F(t)=\int_{\Omega}(u+\xi(v-1-\ln v)) dx.
$$
Then we have
$$
F'(t)+\int_{\Omega}\frac{|\nabla v|^2}{v^2}dx+\mu\int_{\Omega}(v-1)^2 dx=(\xi-\rho)\int_{\Omega} u dx\le 0.
$$
Noting that $F(t)\ge 0$, then we further have
\begin{equation}
\label{4-24}
\int_1^\infty\int_{\Omega}\frac{|\nabla v|^2}{v^2}dxdt+\mu\int_1^\infty\int_{\Omega}(v-1)^2 dxdt\le F(1)\le C.
\end{equation}
Multiplying the second equation of \eqref{1-1} by $v-1$, we obtain
\begin{align*}
\frac{d}{dt}\int_{\Omega}|v-1|^2dx+2\int_{\Omega}|\nabla v|^2 dx+2\mu\int_{\Omega}v(v-1)^2dx\le C,
\end{align*}
which means that
$$
\int_{\Omega}|v(x,t)-1|^2dx-\int_{\Omega}|v(x,s)-1|^2dx\le C(t-s),  \quad \text{for any $1\le s<t$}.
$$
Hence, by Lemma \ref{lem4-1} and  \eqref{4-24}, we have
$$
\lim_{t\to\infty}\|v(\cdot, t)-1\|_{L^2}=0,
$$
and moreover,
$$
\|v(\cdot, t)-1\|_{L^\infty}\le C_1\|v(\cdot, t)-1\|_{L^2}^{\frac25}\|\nabla v\|_{L^\infty}^{\frac35}+
C_2\|v(\cdot, t)-1\|_{L^2}\to 0, \quad\text{as}\quad t\to\infty.
$$
Note that $v$ is H\"older continuous since $v\in W_p^{2,1}$ for $p>5$, which implies that for any
$\varepsilon\in(0,\rho-\xi)$, there exists $T>0$ such that
$$
\xi v(x,t)-\rho<-\varepsilon, \quad \text{for all $t\ge T$ and $x\in\overline\Omega$}.
$$
Then for any $t\ge T$, we have
$$
\frac{d}{dt}\int_\Omega udx=\int_{\omega}(\xi v-\rho)udx\le -\varepsilon\int_\Omega udx,
$$
which implies that
$$
\int_\Omega u(x,t)dx\le e^{-\varepsilon(t-T)}\int_\Omega u(x,T)dx\le Ce^{-\varepsilon t},\qquad\text{for any $t\ge T$},
$$
and we further have
$$
\lim_{t\to\infty}\|u\|_{L^p}=0 \quad \text{for any $p>1$}.
$$
The proof is complete. \hfill $\Box$

By Lemma \ref{lem4-6}, we further have
\begin{lemma}
\label{lem4-7}
Assume that $\mu>0$, $\rho>0$, $0\le \xi< \rho$, and  $v_0\not\equiv 0$. Let $(u, v)$ with $u\in\mathcal X_1, v\in
\mathcal X_2$  be the global solution. Under the assumptions of Lemma \ref{lem4-5} or Lemma \ref{lem4-6}.
We also have
\begin{align}
\label{4-26}
\lim_{t\to\infty}\|u\|_{L^\infty}=0.
\end{align}
\end{lemma}

{\it\bfseries Proof.} We use Moser iteration technique to show \eqref{4-26}.

Taking a cut off function $\eta(t)\in C^1[T, +\infty)$ with $\eta(T)=0$ and  $|\eta'(t)|\le 1$, and $\eta(t)=1$ for $t>T+1$.
Multiplying the first equation of \eqref{1-1} by $r\eta^r u^{r-1}$ for any $r>4m$,
then integrating it over $\Omega$, and using \eqref{1-3}, we obtain
\begin{align}
&\frac{d}{dt}\int_\Omega \eta^ru^{r} dx+mr(r-1)\int_\Omega\eta^r u^{m+r-3}|\nabla u|^2 dx
+\rho r\int_\Omega \eta^{r}u^{r} dx+\int_\Omega \eta^{r}u^{r} dx\nonumber
\\
&= \chi r(r-1)\int_\Omega\eta^ru^{r-1}\nabla v\nabla u dx+
+r\int_\Omega \eta' \eta^{r-1}u^{r} dx+\xi r\int_\Omega \eta^r u^r vdx+\int_\Omega \eta^{r}u^{r} dx\nonumber
\\
\label{4-27}
&\le  \frac{mr(r-1)}4\int_\Omega\eta^r u^{m+r-3}|\nabla u|^2 dx+
Cr^2\int_\Omega\eta^ru^{r+1-m} dx+Cr\int_\Omega\eta^{r-1} u^r dx.
\end{align}
By the boundedness of $u$, using Gagliardo-Nirenberg interpolation inequality,
for any sufficiently small $\sigma>0$, we obtain
\begin{align*}
C r^2\eta^r\|u\|_{L^{r+1-m}}^{r+1-m}
=&C r^2\eta^r\|u^{\frac{r+m-1}2}\|_{L^{\frac{2(r+1-m)}{r+m-1}}}^{\frac{2(r+1-m)}{r+m-1}}
\\
\le &C_1r^2\eta^r\|\nabla u^{\frac{r+m-1}2}\|_{L^2}^{\frac{6(r+2-2m)}{6(m-1)+5r}}
\|u^{\frac{r+m-1}2}\|_{L^{\frac{r}{m+r-1}}}^{\frac{4r(r+2(m-1))}{(6(m-1)+5r)(r+m-1)}}
+C_2r^2\eta^r\|u\|_{\frac r2}^{r+1-m}
\\
\le & \sigma\eta^r\|\nabla u^{\frac{r+m-1}2}\|_{L^2}^2+C_\sigma r^{\frac{6(m-1)+5r}{6(m-1)+r}}
\eta^r\|u\|_{\frac r2}^{\frac{r(2(m-1)+r)}{6(m-1)+r}}+
C_2r^2\eta^r\|u\|_{\frac r2}^{r+1-m}
\\
\le & \sigma\eta^r\|\nabla u^{\frac{r+m-1}2}\|_{L^2}^2+C_{\sigma}r^5\eta^r
\|u\|_{\frac r2}^{\frac{r(2(m-1)+r)}{6(m-1)+r}}+
C_2r^2\eta^r\|u\|_{\frac r2}^{r+1-m}
\\
\le &\sigma\eta^r\|\nabla u^{\frac{r+m-1}2}\|_{L^2}^2+\hat C_{\sigma}r^5
\|\eta u\|_{\frac r2}^{r-4m},
\end{align*}
and
\begin{align*}
Cr\eta^{r-1}\|u\|_{L^r}^{r}=&Cr\eta^{r-1}\|u^{\frac{r+m-1}2}\|_{L^{\frac{2r}{r+m-1}}}^{\frac{2r}{r+m-1}}
\\
\le &C_1r\eta^{r-1}\|\nabla u^{\frac{r+m-1}2}\|_{L^2}^{\frac{6r}{6(m-1)+5r}}
\|u^{\frac{r+m-1}2}\|_{L^{\frac{r}{m+r-1}}}^{\frac{4r^2+6r(m-1)}{(6(m-1)+5r)(r+m-1)}}
+C_2r\eta^{r-1}\|u\|_{\frac r2}^{r}
\\
\le & \sigma \eta^r\|\nabla u^{\frac{r+m-1}2}\|_{L^2}^2+C_\sigma \eta^{r-\frac{6(m-1)+5r}{6(m-1)+2r}}r^{\frac{6(m-1)+5r}{6(m-1)+2r}}
\|u\|_{\frac r2}^{\frac{r(3(m-1)+2r)}{6(m-1)+2r}}+
C_2r\eta^{r-1}\|u\|_{\frac r2}^{r}
\\
\le & \sigma\eta^r\|\nabla u^{\frac{r+m-1}2}\|_{L^2}^2+\tilde C_{\sigma}r^{5/2}
\|\eta u\|_{\frac r2}^{r-3m}+
C_2r\|\eta u\|_{\frac r2}^{r-1}.
\end{align*}
Taking $\sigma$ appropriately small in the above inequalities, and substituting them into \eqref{4-27}, we obtain
\begin{equation}
\label{4-28}
\frac{d}{dt}\int_\Omega \eta^ru^{r} dx+\int_\Omega \eta^{r}u^{r} dx\le C_3r^5\|\eta u\|_{\frac r2}^{r-4m}.
\end{equation}
Letting $r_j=2r_{j-1}=2^j r_0$, $r_0=5m$,
$M_j=\displaystyle\sup_{t\in(T,\infty)}\|\eta u\|_{L^{r_j}}$, and noticing that $\eta(T)=0$,
then by a direct calculation, we conclude that
$$
M_j^{r_j}\le C_3r_j^5M_{j-1}^{r_j-4m}.
$$
By an iteration process, we see that
\begin{align*}
M_j\le & C_3^{\frac 1{r_j}}r_j^{\frac 5{r_j}}M_{j-1}^{1-\frac{4m}{r_02^j}}
\le  C_3^{\sum_{k=1}^j\frac 1{r_02^k}}
r_0^{\sum_{k=1}^j\frac {5}{r_02^k}}M_0^{\prod_{k=1}^j (1-\frac{4m}{r_02^k})}
\le C_4M_0^{\prod_{k=1}^j\frac{r_02^k-4m}{r_02^k}}.
\end{align*}
Next, we show that $S=\prod_{k=1}^\infty\frac{r_02^k-4m}{r_02^k}>0$.
Notice that $\ln \frac1S=\sum_{k=1}^\infty\ln(1+\frac{4m}{r_02^k-4m})$, clearly,
$\sum_{k=1}^\infty\ln(1+\frac{4m}{r_02^k-4m})$ converges, it implies that $S>0$. Then
$$
M_j\le  C_4M_0^{S},\  \text{for any }\ j\ge 1.
$$
Letting $j\to\infty$, we obtain
$$
\sup_{t>T}\|\eta u(x,t)\|_{L^\infty}\le  C_4\sup_{t>T}\|\eta u(x,t)\|_{L^{5m}}^{S},
$$
and by lemma \ref{4-6}, this lemma is proved. \hfill $\Box$

Then Theorem \ref{thm-3} is a direct result of Lemma \ref{lem4-4}--Lemma \ref{lem4-7}.

\section{Extend to the coupled chemotaxis-Stokes system}
In this section, we extend the results to the chemotaxis-Stokes system \eqref{1-13}.
Similar to Section 3, to study the existence of solutions to the system \eqref{1-13}, we also
consider the following approximate problems
\begin{align}
\label{5-1}\left\{
\begin{aligned}
&n_{\varepsilon t}+u_{\varepsilon}\cdot\nabla n_{\varepsilon}=\Delta (\varepsilon n_{\varepsilon} +n_{\varepsilon}^m)-\chi
\nabla\cdot\left(n_{\varepsilon}\nabla c_{\varepsilon}\right)-\varepsilon n_{\varepsilon}^2,  \  (x,t) \in Q,,
\\
&c_{\varepsilon t}+u_{\varepsilon}\cdot\nabla c_{\varepsilon}-\Delta c_{\varepsilon}=-c_{\varepsilon}n_{\varepsilon}, \   (x,t) \in Q,
\\
&u_{\varepsilon t}+\nabla P_{\varepsilon}=\Delta u_{\varepsilon}+n_{\varepsilon}\nabla\varphi_{\varepsilon}, \    (x,t) \in Q,
\\
& {\rm div} u_{\varepsilon}=0, \    (x,t) \in Q,
\\
&\left.\frac{\partial n_{\varepsilon}}{\partial {\bf n}}\right|_{\partial\Omega}
=\left.\frac{\partial c_{\varepsilon}}{\partial {\bf n}}\right|_{\partial\Omega}=0,,
\quad u_{\varepsilon}|_{\partial\Omega}=0,
\\
& n_{\varepsilon}(x,0)=n_{\varepsilon 0}(x)\ge 0,  c_{\varepsilon}(x,0)=c_{\varepsilon 0}(x)\ge 0,
u_{\varepsilon}(x,0)=u_{\varepsilon 0}(x), \quad  x\in\Omega,
\end{aligned}\right.
\end{align}
where $\varphi_\varepsilon\in C^{1+\alpha,\frac\alpha 2}(\overline{\Omega}\times[0,+\infty))$,
$n_{\varepsilon 0}, c_{\varepsilon 0}, u_{\varepsilon 0} \in C^{2+\alpha}(\overline\Omega)$ with
$\left.\frac{\partial n_{\varepsilon 0}}{\partial{\bf n}}\right|_{\partial\Omega}=
0, \left.\frac{\partial c_{\varepsilon 0}}{\partial{\bf n}}\right|_{\partial\Omega}=0$,
$\left.u_{\varepsilon 0}\right|_{\partial\Omega}=0$,
$$
n_{\varepsilon0} \rightarrow n_0, c_{\varepsilon0}\rightarrow c_0,
u_{\varepsilon0}\rightarrow u_0, \nabla\varphi_\varepsilon\rightarrow\nabla\varphi,\quad
\text{uniformly},
$$
\begin{align*}
&\|n_{\varepsilon0}\|_{L^\infty}+\|\nabla n_{\varepsilon0}^m\|_{L^2}
+\|c_{\varepsilon0}\|_{W^{2,\infty}}+\|Au_{\varepsilon0}\|_{L^2}
+\|u_{\varepsilon0}\|_{L^\infty}+\|\nabla\varphi_\varepsilon\|_{L^\infty}
\\
\le &2(\|n_{0}\|_{L^\infty}+\|\nabla n_{0}^m\|_{L^2}
+\|c_{0}\|_{W^{2,\infty}}+\|Au_{0}\|_{L^2}
+\|u_{0}\|_{L^\infty}+\|\nabla\varphi\|_{L^\infty}).
\end{align*}
where $A$ is the Stokes operator, that is $Aw:=-P\Delta w$,
$P:L^r(\Omega)\rightarrow L^r_\sigma(\Omega)$ is the Helmholtz projection  \cite{G1}.
And  $A$ generates a bounded
analytic semigroup $\{e^{-tA}\}_{t\ge 0}$ on $L_\sigma^r$, and the solution $u$ of \eqref{1-1}
can be expressed as
\begin{equation}
\label{5-2}
u=e^{-tA}u_0+\int_{0}^t e^{-(t-s)A}P(n(s)\nabla\varphi(s))ds.
\end{equation}
For more details of Stokes operator, we refer to \cite{KY}.
According to the arguments in \cite{YJ},  each of these problems  is globally
solvable in the classical sense. That is

\begin{lemma}
\label{lem5-1}
Assume that $m>1$, then for any $\varepsilon>0$, the problem \eqref{5-1} admits a unique
classical solution $(n_\varepsilon, c_\varepsilon, u_\varepsilon, \pi_\varepsilon)$
with $n_\varepsilon\ge 0, c_\varepsilon\ge 0$ and $n_\varepsilon, c_\varepsilon,
u_\varepsilon, \in C^{2+\alpha, 1+\frac\alpha 2}(\overline\Omega\times[0,+\infty))$,
$\pi_\varepsilon \in C^{1+\alpha, \frac\alpha 2}(\overline\Omega\times[0,+\infty))$.
\end{lemma}

Completely similar to Lemma \ref{lem3-2}, it is easy to obtain the following lemma.
\begin{lemma}
\label{lem5-2}
Let $(n_\varepsilon, c_\varepsilon, u_\varepsilon, \pi_\varepsilon)$ be the classical
solution of \eqref{5-1}, then we have
\begin{align}
\label{5-3}
&\sup_{t\in(0,\infty)}\|c_\varepsilon(\cdot,t)\|_{L^\infty}\le C_1,
\\
\label{5-4}
&\sup_{t\in(0,\infty)}\int_\Omega n_\varepsilon dx\le C_2,
\end{align}
where $C_1$, $C_2$ are independent of $\varepsilon$.
\end{lemma}

Similar to Lemma \ref{lem3-3}, we also have
\begin{lemma}
\label{lem5-3}
Let $(n_\varepsilon, c_\varepsilon, u_\varepsilon, \pi_\varepsilon)$ be the classical
solution of \eqref{5-1}. Then we have
\begin{align}
&\sup_{t\in(0,+\infty)}\int_\Omega\left(\frac{|\nabla c_\varepsilon|^2}{c_\varepsilon}+
n_{\varepsilon}\ln n_{\varepsilon}+u_{\varepsilon}^2\right)dx+\sup_{t\in(0,+\infty)}\int_{t}^{t+1}
\int_\Omega\left(\varepsilon n_{\varepsilon}^2\ln(1+ n_{\varepsilon})+
\varepsilon\frac{|\nabla n_{\varepsilon}|^2}{n_{\varepsilon}}
\right) dx
\nonumber
\\
\label{5-6}
&+\sup_{t\in(0,+\infty)}\int_{t}^{t+1}
\int_\Omega \left(c_{\varepsilon}|D^2\ln c_{\varepsilon}|^2 dx
+|\nabla n_{\varepsilon}^{\frac m2}|^2+\frac{n_\varepsilon}{c_\varepsilon}|\nabla c_\varepsilon|^2
+\frac{|\nabla c_{\varepsilon}|^4}{c_{\varepsilon}^3}+n_{\varepsilon}^{m+\frac 23}
+|\nabla u_{\varepsilon}|^2\right)dx\le C,
\end{align}
where $C$ is independent of $\varepsilon$.
\end{lemma}

{\it\bfseries Proof.} Similar to \eqref{3-7}, we obtain
\begin{align}
&\frac12\frac{d}{dt}\int_\Omega\frac{|\nabla c_\varepsilon|^2}{c_\varepsilon}dx
+\int_\Omega c_{\varepsilon}|D^2\ln c_{\varepsilon}|^2 dx+
\frac{1}2\int_\Omega
\frac{n_\varepsilon}{c_\varepsilon}|\nabla c_\varepsilon|^2dx
+\int_\Omega
\nabla n_{\varepsilon}\nabla c_{\varepsilon}dx \nonumber
\\
&=\frac12\int_{\partial\Omega}\frac{1}{c_\varepsilon}
\frac{\partial}{\partial{\bf n}}|\nabla c_\varepsilon|^2 ds
+\int_\Omega\frac{u_\varepsilon\cdot\nabla c_\varepsilon}{c_\varepsilon}\Delta c_\varepsilon dx
-\frac12\int_\Omega u_\varepsilon\cdot\nabla c_\varepsilon\frac{|\nabla c_\varepsilon|^2}{c_\varepsilon^2}dx
\nonumber
\\
\label{5-7}
&=\frac12\int_{\partial\Omega}\frac{1}{c_\varepsilon}
\frac{\partial}{\partial{\bf n}}|\nabla c_\varepsilon|^2 ds
-\int_\Omega\frac{\nabla c_\varepsilon\cdot\nabla u_\varepsilon\cdot\nabla c_\varepsilon}{c_\varepsilon}dx.
\end{align}
Noticing that ${\rm div}u_\varepsilon=0$, then for $n_\varepsilon$, multiplying the first equation of \eqref{5-1} by $1+\ln n_{\varepsilon}$, we get
\begin{align}\label{5-8}
\frac{d}{dt}\int_\Omega n_{\varepsilon}\ln n_{\varepsilon}dx+\varepsilon\int_\Omega
\frac{|\nabla n_{\varepsilon}|^2}{n_{\varepsilon}}dx
+\frac{4}{m}\int_\Omega
|\nabla n_{\varepsilon}^{\frac m2}|^2dx
+\varepsilon\int_\Omega n_{\varepsilon}^2(1+\ln n_{\varepsilon})dx
=\chi\int_\Omega \nabla n_{\varepsilon}\nabla c_{\varepsilon} dx.
\end{align}
Comparing with the proof of Lemma \ref{lem3-3}, we only need to deal with the last term of \eqref{5-7}, and we have
\begin{align}
&\frac{d}{dt}\int_\Omega\left(\frac12\frac{|\nabla c_\varepsilon|^2}{c_\varepsilon}+
\frac{1}{\chi}n_{\varepsilon}\ln n_{\varepsilon}\right)dx
+\frac12\int_\Omega c_{\varepsilon}|D^2\ln c_{\varepsilon}|^2 dx+\frac{\varepsilon}{\chi}\int_\Omega
\frac{|\nabla n_{\varepsilon}|^2}{n_{\varepsilon}}dx+\frac{2}{m\chi}\int_\Omega
|\nabla n_{\varepsilon}^{\frac m2}|^2dx\nonumber
\\
&+\frac{\varepsilon}{\chi}\int_\Omega n_{\varepsilon}^2\ln(1+ n_{\varepsilon})dx
+\frac{1}2\int_\Omega \frac{n_\varepsilon}{c_\varepsilon}|\nabla c_\varepsilon|^2dx\nonumber
\\
& \le -\int_\Omega\frac{\nabla c_\varepsilon\cdot\nabla u_\varepsilon\cdot\nabla c_\varepsilon}{c_\varepsilon}dx+
C\nonumber
\\
\label{5-9}
& \le \eta\int_\Omega \frac{|\nabla c_\varepsilon|^4}{c_\varepsilon^3}dx+C_\eta\int_\Omega
|\nabla u_\varepsilon|^2 dx+C
\end{align}
for any sufficiently small $\eta>0$, where $C_\eta$ depends on $\eta$.
By \eqref{2-6}, then we further have
\begin{align}
&\frac{d}{dt}\int_\Omega\left(\frac12\frac{|\nabla c_\varepsilon|^2}{c_\varepsilon}+
\frac{1}{\chi}n_{\varepsilon}\ln n_{\varepsilon}\right)dx
+\frac14\int_\Omega c_{\varepsilon}|D^2\ln c_{\varepsilon}|^2 dx+\frac{\varepsilon}{\chi}\int_\Omega
\frac{|\nabla n_{\varepsilon}|^2}{n_{\varepsilon}}dx\nonumber
\\
&+\frac{2}{m\chi}\int_\Omega
|\nabla n_{\varepsilon}^{\frac m2}|^2dx+\frac{\varepsilon}{\chi}\int_\Omega n_{\varepsilon}^2\ln(1+ n_{\varepsilon})dx
+\int_\Omega
+\frac{1}2\frac{n_\varepsilon}{c_\varepsilon}|\nabla c_\varepsilon|^2 dx\nonumber
\\
\label{5-10}
& \le \hat C\int_\Omega |\nabla u_\varepsilon|^2 dx+C.
\end{align}
Multiplying the third equation of \eqref{5-1} by $u_\varepsilon$, and integrating it over $\Omega$, we obtain
\begin{align}
&\frac12\frac{d}{dt}\int_\Omega u_\varepsilon^2 dx+\int_\Omega |\nabla u_\varepsilon|^2 dx
=\int_\Omega n_\varepsilon\nabla\varphi \cdot u_\varepsilon dx
\le \|\nabla\varphi\|_{L^\infty}\|n_\varepsilon\|_{L^{\frac 65}}\|u_\varepsilon\|_{L^6}\nonumber
\\
&\le C\|n_\varepsilon\|_{L^{\frac 65}}\|\nabla u_\varepsilon\|_{L^2}\nonumber
\\
\label{5-11}
&\le \frac12\int_\Omega |\nabla u_\varepsilon|^2 dx +\tilde C\|n_\varepsilon\|_{L^{\frac 65}}^2.
\end{align}
By Gagliardo-Nirenberg interpolation inequality, we obtain
$$
\tilde C\|n_\varepsilon\|_{L^{\frac 65}}^2=\tilde C\|n_\varepsilon^{\frac m2}\|_{L^{\frac{12}{5m}}}^{\frac 4m}
\le C_1\|\nabla n_\varepsilon^{\frac m2}\|_{L^2}^{\frac{2}{3m-1}}\|n_\varepsilon^{\frac m2}\|_{L^{\frac2m}}+
C_2\|n_\varepsilon\|_{L^1}^2
\le C_3\|\nabla n_\varepsilon^{\frac m2}\|_{L^2}^{\frac{2}{3m-1}}+C_4,
$$
substituting it into \eqref{5-11}, we obtain
\begin{align}\label{5-12}
\frac{d}{dt}\int_\Omega u_\varepsilon^2 dx+\int_\Omega |\nabla u_\varepsilon|^2 dx
\le 2C_3\|\nabla n_\varepsilon^{\frac m2}\|_{L^2}^{\frac{2}{3m-1}}+2C_4.
\end{align}
Combining \eqref{5-10} and \eqref{5-12}, and noticing that ${\frac{2}{3m-1}}<2$, then we obtain
\begin{align*}
&\frac{d}{dt}\int_\Omega\left(\frac12\frac{|\nabla c_\varepsilon|^2}{c_\varepsilon}+
\frac{1}{\chi}n_{\varepsilon}\ln n_{\varepsilon}+2\hat C u_\varepsilon^2\right)dx
+\frac14\int_\Omega c_{\varepsilon}|D^2\ln c_{\varepsilon}|^2 dx+\frac{\varepsilon}{\chi}\int_\Omega
\frac{|\nabla n_{\varepsilon}|^2}{n_{\varepsilon}}dx
\\
&+\frac{2}{m\chi}\int_\Omega
|\nabla n_{\varepsilon}^{\frac m2}|^2dx+\frac{\varepsilon}{\chi}\int_\Omega n_{\varepsilon}^2\ln(1+ n_{\varepsilon})dx
+\int_\Omega \left(\frac{1}2\frac{n_\varepsilon}{c_\varepsilon}|\nabla c_\varepsilon|^2+
\hat C |\nabla u_\varepsilon|^2 \right)dx
\\
& \le 4\hat C C_3\|\nabla n_\varepsilon^{\frac m2}\|_{L^2}^{\frac{2}{3m-1}}+C_5
\\
&\le \frac{1}{m\chi}\|\nabla n_\varepsilon^{\frac m2}\|_{L^2}^2+C_6,
\end{align*}
that is
\begin{align}
&\frac{d}{dt}\int_\Omega\left(\frac12\frac{|\nabla c_\varepsilon|^2}{c_\varepsilon}+
\frac{1}{\chi}n_{\varepsilon}\ln n_{\varepsilon}+2\hat C u_\varepsilon^2\right)dx
+\frac14\int_\Omega c_{\varepsilon}|D^2\ln c_{\varepsilon}|^2 dx+\frac{\varepsilon}{\chi}\int_\Omega
\frac{|\nabla n_{\varepsilon}|^2}{n_{\varepsilon}}dx\nonumber
\\
\label{5-13}
&+\frac{1}{m\chi}\int_\Omega
|\nabla n_{\varepsilon}^{\frac m2}|^2dx +\frac{\varepsilon}{\chi}\int_\Omega n_{\varepsilon}^2\ln(1+ n_{\varepsilon})dx
+\int_\Omega \left(\frac{1}2\frac{n_\varepsilon}{c_\varepsilon}|\nabla c_\varepsilon|^2+
\hat C |\nabla u_\varepsilon|^2 \right)dx \le C.
\end{align}
Recalling \eqref{3-10}, and using Lemma \ref{lem2-2}, we obtain \eqref{5-6}. \hfill $\Box$

Next, We also separate out the case $m>2$.  Similar to Lemma \ref{lem3-4}, we obtain
\begin{lemma}
\label{lem5-4}
Let $(n_\varepsilon, c_\varepsilon, u_\varepsilon, \pi_\varepsilon)$ be the classical
solution of \eqref{5-1}.  If $m>2$,  then we have
\begin{equation}
\label{5-14}
\sup_{t\in(0, +\infty)}\left(\|n_\varepsilon(\cdot,t)\|_{L^\infty}+\|c_\varepsilon(\cdot,t)\|_{W^{1,\infty}} \right)
\le C,
\end{equation}
and for any $r>0$,
\begin{equation}
\label{5-15}
\sup_{t\in(0, +\infty)}\|c_\varepsilon\|_{W_r^{2,1}Q_1(t)}\le C_r,
\end{equation}
where  $C_r$ depends on $r$, both $C$ and $C_r$ are independent of $\varepsilon$.
\end{lemma}

{\it\bfseries Proof. } Recalling \eqref{5-6}, we see that
\begin{equation}
\label{5-16}
\sup_{t\in(0,infty)}\int_t^{t+1}\int_\Omega n_\varepsilon^{m+\frac23} dxds\le C.
\end{equation}
Multiplying the third equation of \eqref{5-1} by $\Delta u_\varepsilon$, and integrating it over $\Omega$, we obtain
\begin{align*}
&\frac12\frac{d}{dt}\int_\Omega |\nabla u_\varepsilon|^2 dx+\int_\Omega |\Delta u_\varepsilon|^2 dx
=\int_\Omega n_\varepsilon\nabla\varphi \cdot \Delta u_\varepsilon dx
\le \|\nabla\varphi\|_{L^\infty}\|n_\varepsilon\|_{L^{2}}\|\Delta u_\varepsilon\|_{L^2}
\\
&\le C\|n_\varepsilon\|_{L^2}^2+\frac12\|\Delta u_\varepsilon\|_{L^2}^2.
\end{align*}
By Lemma \ref{lem2-2}, and using \eqref{5-16}, and noticing that $m>2$,
then we obtain
\begin{equation}
\label{5-17}
\sup_{t\in(0,\infty)}\|\nabla u_\varepsilon(\cdot, t)\|_{L^2}^2+
\sup_{t\in(0,\infty)}\int_t^{t+1}\int_\Omega|\Delta u_\varepsilon|^2 dxds
\le C.
\end{equation}
Completely similar to \eqref{3-20}, for any $r>0$, we have
\begin{align}
&\frac1{r+1}\frac{d}{dt}\int_\Omega n_\varepsilon^{r+1} dx
+\frac{rm}2\int_\Omega n_\varepsilon^{m+r-2}|\nabla n_\varepsilon|^2 dx+\frac12\int_\Omega n_\varepsilon^{r+1} dx
\nonumber
\\
\label{5-18}
\le & C\int_\Omega|\Delta c_\varepsilon|^{r+3-m}dx+C.
\end{align}
By Lemma \ref{lem2-2}, we get
\begin{align}\label{5-19}
\sup_{t\in(0, +\infty)}\|n_\varepsilon(\cdot,t)\|_{L^{r+1}}^{r+1} +
\sup_{t\in(0, +\infty)}\int_t^{t+1}\|\nabla n_\varepsilon^{\frac{m+r}2}(\cdot,s)\|_{L^{2}}^{2} ds
\le C\sup_{t\in(0, +\infty)}\int_t^{t+1}\int_\Omega|\Delta c_\varepsilon|^{r+3-m}dxds+C.
\end{align}
By the second equation of \eqref{5-1}, we see that
\begin{align}
\label{5-20}
c_{\varepsilon t}-\Delta c_{\varepsilon}+c_{\varepsilon}=-u_{\varepsilon}\nabla c_{\varepsilon}+
c_{\varepsilon}-c_{\varepsilon}n_{\varepsilon}.
\end{align}
By Lemma \ref{lem2-3}, and using \eqref{5-3}, \eqref{5-17}
and Gagliardo-Nirenberg interpolation inequality, we have
\begin{align*}
&\sup_{t\in(0, +\infty)}\int_t^{t+1}\int_\Omega|\Delta c_\varepsilon|^{3}dxds
\le C+C\sup_{t\in(0, +\infty)}\int_t^{t+1}(\|u_{\varepsilon}\cdot\nabla c_{\varepsilon}\|_{L^{3}}^{3}
+\|n_{\varepsilon}\|_{L^{3}}^{3}) ds
\\
&\le C+C\sup_{t\in(0, +\infty)}\int_t^{t+1}(\|u_{\varepsilon}\|_{L^6}^3\|\nabla c_{\varepsilon}\|_{L^6}^3
+\|n_{\varepsilon}\|_{L^{3}}^{3}) ds
\\
&\le C+\tilde C\sup_{t\in(0, +\infty)}\int_t^{t+1}(1+\|\Delta c_{\varepsilon}\|_{L^3}^{3/2}
+\|n_{\varepsilon}\|_{L^{3}}^{3}) ds,
\end{align*}
which implies that
\begin{align}
\label{5-21}
\sup_{t\in(0, +\infty)}\int_t^{t+1}\int_\Omega|\Delta c_\varepsilon|^{3}dxds
\le C_1+C_2\sup_{t\in(0, +\infty)}\int_t^{t+1}\|n_{\varepsilon}\|_{L^{3}}^{3} ds.
\end{align}
Substituting this inequality into \eqref{5-19} with $r=m$, then
we obtain
\begin{align}\label{5-22}
\sup_{t\in(0, +\infty)}\|n_\varepsilon(\cdot,t)\|_{L^{m+1}}^{m+1} +
\sup_{t\in(0, +\infty)}\int_t^{t+1}\|\nabla n_\varepsilon^m(\cdot,s)\|_{L^{2}}^{2} ds
\le C\sup_{t\in(0, +\infty)}\int_t^{t+1}\|n_{\varepsilon}\|_{L^{3}}^{3}ds+C.
\end{align}
Noticing that $m>2$, it implies that
\begin{align}\label{5-23}
\sup_{t\in(0, +\infty)}\|n_\varepsilon(\cdot,t)\|_{L^{m+1}}^{m+1} +
\sup_{t\in(0, +\infty)}\int_t^{t+1}\|\nabla n_\varepsilon^m(\cdot,s)\|_{L^{2}}^{2} ds
\le C.
\end{align}
By Duhamel's principle,
we see that the solution of \eqref{5-20} can be expressed as follows
$$
c_\varepsilon(x,t)=e^{-t}e^{t\Delta}c_{\varepsilon 0}+\int_{0}^t e^{-(t-s)}e^{(t-s)\Delta}(
-u_{\varepsilon}\nabla c_{\varepsilon}+
c_{\varepsilon}-c_{\varepsilon}n_{\varepsilon})ds,
$$
where $\{e^{t\Delta}\}_{t\ge 0}$ is the heat semigroup on the domain $\Omega$ under Neumann boundary condition,
for more properties of
Neumann heat semigroup, please refer to \cite{W4}.
Then by \eqref{5-3}, \eqref{5-6}, \eqref{5-17}, \eqref{5-22}, and noticing that $m>2$, we have
\begin{align*}
&\|\nabla c_\varepsilon(\cdot, t)\|_{L^\infty}\le  e^{-t}\|\nabla c_{\varepsilon 0}\|_{L^\infty}
+\int_0^t e^{-(t-s)}(t-s)^{-\frac 32(\frac 1{m+1})-\frac12}\|c_{\varepsilon}n_{\varepsilon}\|_{L^{m+1}}ds
\\
&+\int_0^t e^{-(t-s)}(t-s)^{-\frac 32(\frac 14)-\frac12}\left(
\|u_{\varepsilon}\cdot\nabla c_{\varepsilon}\|_{L^4}+\|c_{\varepsilon}\|_{L^4}\right)ds
\\
\le &  e^{-t}\|\nabla c_{\varepsilon 0}\|_{L^\infty}
+C_1\int_0^t e^{-(t-s)}(t-s)^{-\frac 32(\frac 1{m+1})-\frac12} ds
+\int_0^t e^{-(t-s)}(t-s)^{-\frac 78}\left(
\|u_{\varepsilon}\|_{L^6}\|\nabla c_{\varepsilon}\|_{L^{12}}+C\right)ds
\\
\le &  e^{-t}\|\nabla c_{\varepsilon 0}\|_{L^\infty}
+C_1\int_0^t e^{-s}s^{-\frac 32(\frac 1{m+1})-\frac12} ds
+\int_0^t e^{-(t-s)}(t-s)^{-\frac 78}\left(C_2
\|\nabla u_{\varepsilon}\|_{L^2}\|\nabla c_{\varepsilon}\|_{L^{2}}^{\frac16}\|\nabla c_{\varepsilon}
\|_{L^{\infty}}^{\frac56}+C\right)ds
\\
\le &  e^{-t}\|\nabla c_{\varepsilon 0}\|_{L^\infty}+
C_1\int_0^\infty e^{-s}s^{-\frac 32(\frac 1{m+1})-\frac12} ds+C_3\sup_{t\in(0, \infty)}\left(
1+\|\nabla c_{\varepsilon}\|_{L^\infty}^{\frac56}\right)
\int_0^{\infty} e^{-s}s^{-\frac 78} ds
\\
\le & e^{-t}\|\nabla c_{\varepsilon 0}\|_{L^\infty}+C_4\left(1+\sup_{t\in(0, T_{\max})}\|\nabla c_{\varepsilon}\|_{L^\infty}^{\frac56}
\right),
\end{align*}
which implies that
\begin{align}\label{5-24}
\|\nabla c_\varepsilon(\cdot, t)\|_{L^\infty}\le C.
\end{align}
Then by Moser iteration method, see for example \cite{J2}, we obtain \eqref{5-14}. \hfill $\Box$

Completely similar to \eqref{3-19}, for any $r>0$, we also have
\begin{align}
&\frac1{r+1}\frac{d}{dt}\int_\Omega n_\varepsilon^{r+1} dx
+\frac{rm}2\int_\Omega n_\varepsilon^{m+r-2}|\nabla n_\varepsilon|^2 dx+\int_\Omega n_\varepsilon^{r+1} dx
\nonumber
\\
\label{5-25}
\le & C\int_\Omega n_\varepsilon^{r+2-m}|\nabla c_\varepsilon|^{2}dx+C.
\end{align}
Taking $r=m-1$ in \eqref{5-25}, and noticing
that $n_\varepsilon|\nabla c_\varepsilon|^{2}\le \|c_\varepsilon\|_{L^\infty}
\frac{n_\varepsilon}{c_\varepsilon}|\nabla c_\varepsilon|^{2}$,
then combining with  \eqref{5-6} yields
\begin{lemma}
\label{lem5-5}
Assume that $m>1$.
Let $(n_\varepsilon, c_\varepsilon, u_\varepsilon, \pi_\varepsilon)$ be the classical
solution of \eqref{5-1}.  Then we have
\begin{equation}
\label{add5-25}
\sup_{t\in(0,+\infty)}\|n_\varepsilon\|_{L^{m}}+\sup_{t\in(0,+\infty)}\int_t^{t+1}
\int_\Omega n_\varepsilon^{2m-3}|\nabla n_\varepsilon|^2 dxds\le C,
\end{equation}
where $C$ is independent of $\varepsilon$.
\end{lemma}
Similar to Lemma \ref{lem3-6}, we obtain

\begin{lemma}
\label{lem5-6}
Assume  $1<m\le 2$.
Let $(n_\varepsilon, c_\varepsilon, u_\varepsilon, \pi_\varepsilon)$ be the classical
solution of \eqref{5-1}. If
$$
\sup_{t\in(0,+\infty)}\int_t^{t+1}
\int_\Omega n_\varepsilon^{\alpha}|\nabla n_\varepsilon|^2 dxds\le C
$$
with $\alpha\ge 0$, then for any $r>0$,
\begin{equation}
\label{5-26}
\sup_{t\in(0, +\infty)}\|n_\varepsilon(\cdot,t)\|_{L^{r+1}}^{r+1} +
\sup_{t\in(0, +\infty)}\int_t^{t+1}\|\nabla n_\varepsilon^{\frac{m+r}2}(\cdot,s)\|_{L^{2}}^{2} ds
\le C_r,
\end{equation}
where $C_r$ depends on $r$, and  which is independent of $\varepsilon$.
\end{lemma}

{\it\bfseries Proof. } Similar to\eqref{3-23}, by \eqref{5-6}, we also have
\begin{equation}
\label{5-27}
\sup_{t\in(0,+\infty)}\int_t^{t+1}
\int_\Omega |\nabla n_\varepsilon|^2 dxds\le C.
\end{equation}
Applying $\nabla$ to the second equation of \eqref{3-1}, and
multiplying the resulting equation by $|\nabla v_\varepsilon|^{r-2}\nabla v_\varepsilon$
for any $r>2$, and using Lemma \ref{lem2-5} and the boundary trace embedding inequalities,
see for example the proof for \eqref{3-25}, we obtain
\begin{align*}
&\frac1r\frac{d}{dt}\int_\Omega |\nabla c_\varepsilon|^{r} dx
+\int_{\Omega}|\nabla c_\varepsilon|^{r-2}|\nabla^2 c_\varepsilon|^2 dx+(r-2)
\int_{\Omega}|\nabla c_\varepsilon|^{r-2}(\nabla|\nabla c_\varepsilon|)^2 dx
+\int_{\Omega} n_\varepsilon|\nabla c_\varepsilon|^{r}dx
\\
&=-\int_{\Omega}c_\varepsilon |\nabla c_\varepsilon|^{r-2}\nabla c_\varepsilon\nabla n_\varepsilon dx
+\int_{\Omega}u_\varepsilon\nabla c_\varepsilon {\rm div}(|\nabla c_\varepsilon|^{r-2}\nabla c_\varepsilon)dx+
\frac12\int_{\partial\Omega}\frac{\partial(|\nabla c_\varepsilon|^2)}{\partial{\bf n}}|\nabla c_\varepsilon|^{r-2} ds
\\
&\le -\int_{\Omega}c_\varepsilon |\nabla c_\varepsilon|^{r-2}\nabla c_\varepsilon\nabla n_\varepsilon dx
+\int_{\Omega}|u_\varepsilon|^2|\nabla c_\varepsilon|^r dx
+\frac12\int_{\Omega}|\nabla c_\varepsilon|^{r-2}|\nabla^2 c_\varepsilon|^2 dx
\\
&+\frac{r-2}4
\int_{\Omega}|\nabla c_\varepsilon|^{r-2}(\nabla|\nabla c_\varepsilon|)^2 dx+
\kappa\int_{\partial\Omega}|\nabla c_\varepsilon|^r ds
\\
&\le -\int_{\Omega}c_\varepsilon |\nabla c_\varepsilon|^{r-2}\nabla c_\varepsilon\nabla n_\varepsilon dx
+\int_{\Omega}|u_\varepsilon|^2|\nabla c_\varepsilon|^r dx
+\frac12\int_{\Omega}|\nabla c_\varepsilon|^{r-2}|\nabla^2 c_\varepsilon|^2 dx
\\
&+\frac{r-2}2
\int_{\Omega}|\nabla c_\varepsilon|^{r-2}(\nabla|\nabla c_\varepsilon|)^2 dx+C,
\end{align*}
that is
\begin{align}
&\frac1r\frac{d}{dt}\int_\Omega |\nabla c_\varepsilon|^{r} dx
+\frac12\int_{\Omega}|\nabla c_\varepsilon|^{r-2}|\nabla^2 c_\varepsilon|^2 dx+\frac{r-2}2
\int_{\Omega}|\nabla c_\varepsilon|^{r-2}(\nabla|\nabla c_\varepsilon|)^2 dx
+\int_{\Omega} n_\varepsilon|\nabla c_\varepsilon|^{r}dx \nonumber
\\
\label{5-28}
&\le -\int_{\Omega}c_\varepsilon |\nabla c_\varepsilon|^{r-2}\nabla c_\varepsilon\nabla n_\varepsilon dx
+\int_{\Omega}|u_\varepsilon|^2|\nabla c_\varepsilon|^r dx+C.
\end{align}
Similar to \eqref{3-26}, we have
\begin{equation}
\label{5-29}
\|\nabla c_\varepsilon\|_{r+2}^{r+2}\le
C\int_\Omega |\nabla c_\varepsilon|^{r-2}|\nabla^2 c_\varepsilon|^2 dx
+C\int_{\Omega}|\nabla c_\varepsilon|^{r-2}(\nabla|\nabla c_\varepsilon|)^2 dx.
\end{equation}
By \eqref{5-4} and \eqref{5-27}, we also have
\begin{equation}
\label{5-30}
\sup_{t\in(0,+\infty)}\int_t^{t+1}
\int_\Omega n_\varepsilon^2 dxds\le C.
\end{equation}
Multiplying the third equation of \eqref{5-1} by $\Delta u_\varepsilon$, and integrating it over $\Omega$,
we obtain
\begin{align*}
\frac1{2}\frac{d}{dt}\int_\Omega|\nabla u_\varepsilon|^2 dx+\int_\Omega|\Delta u_\varepsilon|^2 dx
=&-\int_\Omega n_{\varepsilon}\nabla\varphi_{\varepsilon}\Delta u_\varepsilon dx
\\
\le & \frac12\int_\Omega|\Delta u_\varepsilon|^2 dx+C\int_\Omega|n_\varepsilon|^2 dx,
\end{align*}
which implies that
$$
\frac{d}{dt}\int_\Omega|\nabla u_\varepsilon|^2 dx+\int_\Omega|\Delta u_\varepsilon|^2 dx
+\int_\Omega|\nabla u_\varepsilon|^2 dx\le C\int_\Omega|n_\varepsilon|^2 dx+\int_\Omega|\nabla u_\varepsilon|^2 dx.
$$
Using Lemma \ref{lem2-2}, \eqref{5-6} and \eqref{5-30}, we obtain
\begin{equation}
\label{5-31}
\sup_{t\in(0,+\infty)}\int_\Omega|\nabla u_\varepsilon|^2 dx+
\sup_{t\in(0,+\infty)}\int_t^{t+1}
\int_\Omega |\Delta u_\varepsilon|^2  dxds\le C.
\end{equation}
Taking $r=4$ in \eqref{5-28}, we obtain
\begin{align*}
&\frac14\frac{d}{dt}\int_\Omega |\nabla c_\varepsilon|^{4} dx
+\frac12\int_{\Omega}|\nabla c_\varepsilon|^2|\nabla^2 c_\varepsilon|^2 dx+
\int_{\Omega}|\nabla c_\varepsilon|^2(\nabla|\nabla c_\varepsilon|)^2 dx
+\int_{\Omega} n_\varepsilon|\nabla c_\varepsilon|^{4}dx
\\
&\le -\int_{\Omega}c_\varepsilon |\nabla c_\varepsilon|^2\nabla c_\varepsilon\nabla n_\varepsilon dx
+\int_{\Omega}|u_\varepsilon|^2|\nabla c_\varepsilon|^4 dx+C
\\
&\le \eta\int_{\Omega} |\nabla c_\varepsilon|^6 dx+C_\eta\int_{\Omega}
(|\nabla n_\varepsilon|^2 +|u_\varepsilon|^6) dx +C
\\
&\le \eta\int_{\Omega} |\nabla c_\varepsilon|^6 dx+C_\eta\int_{\Omega}
|\nabla n_\varepsilon|^2 dx +\tilde C_\eta \|\nabla u_\varepsilon\|_{L^2}^6 +C,
\end{align*}
Using \eqref{5-29} and \eqref{5-31}, and taking $\eta$ appropriately small, then we have
\begin{align*}
&\frac14\frac{d}{dt}\int_\Omega |\nabla c_\varepsilon|^{4} dx
+\frac14\int_{\Omega}|\nabla c_\varepsilon|^2|\nabla^2 c_\varepsilon|^2 dx+
\frac14\int_{\Omega}|\nabla c_\varepsilon|^2(\nabla|\nabla c_\varepsilon|)^2 dx
+\int_{\Omega} n_\varepsilon|\nabla c_\varepsilon|^{4}dx+
\sigma\int_\Omega |\nabla c_\varepsilon|^{6} dx
\\
&\le C\int_{\Omega}|\nabla n_\varepsilon|^2 dx +C.
\end{align*}
Using Lemma \ref{lem2-2} and \eqref{5-27}, then we obtain
\begin{align}\label{5-32}
\sup_{t\in(0,+\infty)}\int_\Omega |\nabla c_\varepsilon|^{4} dx\le C.
\end{align}
By using this inequality, and completely similar to the proof of Lemma \ref{lem3-6},
we obtain \eqref{5-26}.   \hfill $\Box$

\begin{lemma}
\label{lem5-7}
Assume  $1<m\le 2$. Let $(n_\varepsilon, c_\varepsilon, u_\varepsilon, \pi_\varepsilon)$ be the classical
solution of \eqref{5-1}. If
\begin{equation}\label{5-33}
\sup_{t\in(0,+\infty)}
\int_\Omega n_\varepsilon^{A_k(m-1)+1} dx+
\sup_{t\in(0,+\infty)}\int_t^{t+1}
\int_\Omega n_\varepsilon^{A_k(m-1)+m-2}|\nabla n_\varepsilon|^2 dxds\le C
\end{equation}
with $A_k(m-1)+m-2< 0$,
then there exists a constant $\tilde C$ independent of $\varepsilon$, such that
\begin{equation}\label{5-34}
\sup_{t\in(0,+\infty)}
\int_\Omega n_\varepsilon^{A_{k+1}(m-1)+1} dx+
\sup_{t\in(0,+\infty)}\int_t^{t+1}
\int_\Omega n_\varepsilon^{A_{k+1}(m-1)+m-2}|\nabla n_\varepsilon|^2 dxds\le \tilde C
\end{equation}
where $A_{k+1}=\frac23(m-1)A_k^2+(\frac{8m}{3}-2)A_k+2m-\frac13$.
\end{lemma}

{\it\bfseries Proof. } Recalling \eqref{5-28}, and  noticing that $2-m-A_k(m-1)>0$,
we see that
\begin{align*}
&\frac1r\frac{d}{dt}\int_\Omega |\nabla c_\varepsilon|^{r} dx
+\frac12\int_{\Omega}|\nabla c_\varepsilon|^{r-2}|\nabla^2 c_\varepsilon|^2 dx+\frac{r-2}2
\int_{\Omega}|\nabla c_\varepsilon|^{r-2}(\nabla|\nabla c_\varepsilon|)^2 dx
+\int_{\Omega} n_\varepsilon|\nabla c_\varepsilon|^{r}dx
\\
&\le -\int_{\Omega}c_\varepsilon |\nabla c_\varepsilon|^{r-2}\nabla c_\varepsilon\nabla n_\varepsilon dx
+\int_{\Omega}|u_\varepsilon|^2|\nabla c_\varepsilon|^r dx+C
\\
&\le-\int_{\Omega}c_\varepsilon |\nabla c_\varepsilon|^{r-2}\nabla c_\varepsilon\nabla n_\varepsilon dx
+\eta\int_\Omega|\nabla c_\varepsilon|^{r+2}dx+C_\eta\int_{\Omega} u_\varepsilon^{r+2}dx
+C=I_1+I_2+I_3+C.
\end{align*}
For $I_1$, it is completely similar to the proof of Lemma \ref{lem3-7};
for $I_2$, using \eqref{5-29},
and taking $\eta$ appropriately small, then it can be controlled by the second
and the third terms on the left side of the above inequality.
Thus, we only need to deal with $I_3$.
Therefore, taking $r=2+2(A_k+1)(m-1)$ in the above inequality, we obtain
\begin{align}
&\frac1{2+2(A_k+1)(m-1)}\frac{d}{dt}\int_\Omega |\nabla c_\varepsilon|^{2+2(A_k+1)(m-1)} dx
+\frac12\int_{\Omega}|\nabla c_\varepsilon|^{2(A_k+1)(m-1)}|\nabla^2 c_\varepsilon|^2 dx\nonumber
\\
&+\frac{(A_k+1)(m-1)}2
\int_{\Omega}|\nabla c_\varepsilon|^{2(A_k+1)(m-1)}(\nabla|\nabla c_\varepsilon|)^2 dx
+\frac12\int_{\Omega} n_\varepsilon|\nabla c_\varepsilon|^{2+2(A_k+1)(m-1)}dx\nonumber
\\
\label{5-35}
&\le  C\int_\Omega n_\varepsilon^{A_k(m-1)+m-2}|\nabla n_\varepsilon|^2 dx
+C\int_\Omega u_\varepsilon^{4+2(A_k+1)(m-1)}dx +C.
\end{align}
For $u_\varepsilon$,  by \eqref{add5-25}, then
\begin{align}
\|u_\varepsilon(\cdot, t)\|_{L^3}\le &e^{-t}\|u_{\varepsilon 0}\|_{L^3}
+\int_0^t \|e^{-(t-s)A}P(n_\varepsilon(s)\nabla\varphi_\varepsilon(s))\|_{L^3}ds\nonumber
\\
\le &e^{-t}\|u_{\varepsilon 0}\|_{L^3}+
\int_0^t e^{-\lambda(t-s)}(t-s)^{-\frac32(\frac1m-\frac13)}\|n_\varepsilon(s)\nabla\varphi_\varepsilon\|_{L^m}ds
\nonumber
\\
\le &e^{-t}\|u_{\varepsilon 0}\|_{L^3}+\sup_{t\in(0,\infty)}
\|n_\varepsilon\|_{L^m}\|\nabla\varphi_\varepsilon\|_{L^\infty}
\int_0^t e^{-\lambda(t-s)}(t-s)^{-\frac32(\frac1m-\frac13)}ds\nonumber
\\
\label{5-36}
\le & C
\end{align}
since that $\frac32(\frac1m-\frac13)<1$.
By \eqref{3-32}, \eqref{5-33}, we also have
\begin{align}
\label{5-37}
\sup_{t\in(0,\infty)}\int_t^{t+1}\int_\Omega n_\varepsilon^{m+\frac23+\frac53 A_k(m-1)}dx ds\le C.
\end{align}
By the $L^{p, q}$ theory of Stokes operator \cite{FK, HP}, we also have
\begin{align}
\label{5-38}
\sup_{t\in(0,\infty)}\int_t^{t+1}
\|A u_\varepsilon\|_{L^{m+\frac23+\frac53 A_k(m-1)}}^{m+\frac23+\frac53 A_k(m-1)}ds\le C.
\end{align}
By Gagliardo-Nirenberg interpolation inequality, we have
$$
\|u_\varepsilon\|_{L^{3m+2+5A_k(m-1)}}^{3m+2+5A_k(m-1)}\le C_1\|u_\varepsilon\|_{L^3}^{2m+\frac43+\frac{10A_k(m-1)}3}
\|A u_\varepsilon\|_{L^{m+\frac23+\frac53 A_k(m-1)}}^{m+\frac23+\frac53 A_k(m-1)}+C_2\|u_\varepsilon\|_{L^3}
^{3m+2+5A_k(m-1)}.
$$
Combining with \eqref{5-36}, \eqref{5-38}, we obtain
\begin{align}
\label{5-39}
\sup_{t\in(0,\infty)}\int_t^{t+1}\|u_\varepsilon\|_{L^{3m+2+5A_k(m-1)}}^{3m+2+5A_k(m-1)}ds\le C.
\end{align}
Recalling \eqref{5-35}, combining with \eqref{5-29}, \eqref{5-33}, \eqref{5-39}, and
noticing that $3m+2+5A_k(m-1)>4+2(A_k+1)(m-1)$, then using Lemma \ref{lem2-2}, we obtain
\begin{align}
&\sup_{t\in(0,\infty)}\int_\Omega |\nabla c_\varepsilon|^{2+2(A_k+1)(m-1)} dx
+\sup_{t\in(0,\infty)}\int_t^{t+1}\int_{\Omega}|\nabla c_\varepsilon|^{2(A_k+1)(m-1)}(|\nabla^2 c_\varepsilon|^2
+(\nabla|\nabla c_\varepsilon|)^2) dxds
\nonumber
\\
\label{5-40}
&+\sup_{t\in(0,\infty)}\int_t^{t+1}\int_{\Omega} n_\varepsilon|\nabla c_\varepsilon|^{2+2(A_k+1)(m-1)}dxds
\le  C.
\end{align}
Recalling \eqref{5-25}, then the following proof is completely similar to \eqref{3-31}-\eqref{3-33},
and the proof is complete. \hfill $\Box$

\begin{lemma}
\label{lem5-8}
Assume  $m>\frac{11}4-\sqrt 3$.
Let $(n_\varepsilon, c_\varepsilon, u_\varepsilon, \pi_\varepsilon)$ be the classical
solution of \eqref{5-1}. Then  we have
\begin{align}
\label{5-41}
&\sup_{t\in(0, +\infty)}\|n_\varepsilon(\cdot,t)\|_{L^{r+1}}^{r+1}+
\sup_{t\in(0, +\infty)}\int_t^{t+1}\|\nabla n_\varepsilon^{\frac{m+r}2}(\cdot,s)\|_{L^{2}}^{2} ds
\le C_r, \ \text{for any $r>0$}
\\
\label{5-42}
&\sup_{t\in(0, +\infty)}\|u_\varepsilon(\cdot,t)\|_{H^1}^{2} +
\sup_{t\in(0, +\infty)}\int_t^{t+1}\|u_\varepsilon(\cdot,s)\|_{H^{2}}^{2} ds
\le C.
\end{align}
where $C_r$, $C$ are independent of $\varepsilon$, $C_r$ depends on $r$.
Using the above two inequalities, we further have
\begin{align}
\label{5-43}
&\sup_{t\in(0,\infty)}(\|u_\varepsilon(\cdot, t)\|_{L^\infty}+\|A^{\beta}u(\cdot, t)\|_{L^2})\le C,
\\
\label{5-44}
&\sup_{t\in(0,\infty)}(\|c_\varepsilon(\cdot, t)\|_{W^{1,\infty}}+\|n_\varepsilon(\cdot, t)\|_{L^\infty})\le C,
\\
\label{5-45}
&\sup_{t\in(0, +\infty)}\|v_\varepsilon\|_{W_p^{2,1}Q_1(t)}\le \tilde C_p, \ \text{for any $p>1$}.
\end{align}
where $C$, $C_p$ are independent of $\varepsilon$, $C_p$ depends on $p$.
\end{lemma}

{\it\bfseries Proof.}
By Lemma \ref{lem3-8}, we see that if $m>\frac{11}4-\sqrt 3$, then $A_k\nearrow +\infty$ as $k\to\infty$.
Using  Lemma \ref{lem5-5}-Lemma \ref{lem5-7}, \eqref{5-41} is arrived.

For $u_\varepsilon$, recalling Lemma \ref{lem5-3} and the proof of \eqref{5-31}, and \eqref{5-42}
is readily obtained. Next, for any $\beta\in(\frac34, 1)$, we note that
\begin{align*}
\|A^{\beta}u_\varepsilon\|_{L^2}\le &e^{-t}\|A^{\beta}u_{\varepsilon 0}\|_{L^2}
+\int_0^t \|A^{\beta}e^{-(t-s)A}P(n_\varepsilon(s)\nabla\varphi_\varepsilon(s))\|_{L^2}ds
\\
\le &e^{-t}\|A^{\beta}u_{\varepsilon 0}\|_{L^2}+\int_0^t e^{-\lambda(t-s)}(t-s)^{-\beta}\|n_\varepsilon(s)
\nabla\varphi_\varepsilon\|_{L^2}ds
\\
\le &e^{-t}\|A^{\beta}u_{\varepsilon 0}\|_{L^2}+\int_0^t e^{-\lambda(t-s)}(t-s)^{-\beta}
\|n_\varepsilon\|_{L^2}\|\nabla\varphi_\varepsilon\|_{L^\infty}ds
\\
\le &C,
\end{align*}
by embedding theorem, $u_\varepsilon\in L^\infty(\Omega\times (0, \infty))$ since $\beta>\frac34$,
and \eqref{5-43} is proved.

For $c_\varepsilon$,  we have
\begin{align*}
\|\nabla c_\varepsilon\|_{L^\infty}\le & e^{-t}\|\nabla c_{\varepsilon 0}\|_{L^\infty}+\int_0^t e^{-(t-s)}
\left\|\nabla\left(e^{(t-s)\Delta}\Big(-u_\varepsilon\cdot\nabla n_\varepsilon-c_\varepsilon n_\varepsilon)\Big)\right)\right\|_{L^\infty}ds
\\
\le & e^{-t}\|\nabla c_{\varepsilon 0}\|_{L^\infty}+\int_0^t e^{-(t-s)}(t-s)^{-\frac12-\frac14}
\left\|-u_\varepsilon\cdot\nabla c_\varepsilon-c_\varepsilon n_\varepsilon\right\|_{L^6}ds
\\
\le & e^{-t}\|\nabla c_{\varepsilon 0}\|_{L^\infty}+\int_0^t e^{-(t-s)}(t-s)^{-\frac34}
(\|u_\varepsilon \|_{L^\infty}\|\nabla c_\varepsilon\|_{L^6}+\|c_\varepsilon\|_{L^\infty}\|n_\varepsilon\|_{L^6})ds
\\
\le & e^{-t}\|\nabla c_{\varepsilon 0}\|_{L^\infty}+\int_0^t e^{-(t-s)}(t-s)^{-\frac34}
(\|u_\varepsilon \|_{L^\infty}\|\nabla c_\varepsilon\|_{L^\infty}^{\frac23}\|\nabla c_\varepsilon\|_{L^\infty}^{\frac13}
+\|c_\varepsilon\|_{L^\infty}\|n_\varepsilon\|_{L^6})ds
\\
\le & e^{-t}\|\nabla c_{\varepsilon 0}\|_{L^\infty}+C(1+\sup_{t\in(0,\infty)}
\|\nabla c_\varepsilon\|_{L^\infty}^{\frac23}).
\end{align*}
Combining with \eqref{5-3}, we obtain
$$
\sup_{t\in(0,\infty)}\|c_\varepsilon(\cdot, t)\|_{W^{1,\infty}}\le C.
$$
Next, by a standard Moser iteration technique, we obtain \eqref{5-44}.
By Lemma \ref{lem2-3}, and \eqref{5-45} is readily arrived.  \hfill $\Box$

By \eqref{5-41}-\eqref{5-42}, we see that $u_\varepsilon\cdot\nabla c_\varepsilon$ is bounded uniformly, then
completely similar to the proof of Lemma \ref{lem3-11}. We also have
\begin{lemma}
\label{lem5-9}
Assume  $m>\frac{11}4-\sqrt 3$.
Let $(n_\varepsilon, c_\varepsilon, u_\varepsilon, \pi_\varepsilon)$ be the classical
solution of \eqref{5-1}. Then  we have
\begin{align}
\label{5-46}
\sup_{t\in(0, +\infty)}\int_\Omega|\nabla n_\varepsilon^m|^2 dx+
\varepsilon\sup_{t\in(0, +\infty)}\int_t^{t+1}\int_\Omega \left|\frac{\partial n_\varepsilon}{\partial t}\right|^2 dxds
+\sup_{t\in(0, +\infty)}\int_t^{t+1}\int_\Omega n_\varepsilon^{m-1}
\left|\frac{\partial n_\varepsilon}{\partial t}\right|^2 dxds\le C,
\end{align}
where $C$ is independent of $\varepsilon$.
\end{lemma}

Similar to the proof of Theorem \ref{thm-1},
by  Lemma \ref{lem5-8}-Lemma \ref{lem5-9}, letting $\varepsilon\to 0$,
and  Theorem \ref{thm-4} is proved.

\end{document}